\documentclass[final,3p]{elsarticle}
 \usepackage{graphics}
 \usepackage{graphicx}
 \usepackage{epsfig}
\usepackage{amssymb}
 \usepackage{amsthm}
 \usepackage{lineno}
 \usepackage{amsmath}
   \numberwithin{equation}{section}
\usepackage{mathrsfs}

\newtheorem{thm}{Theorem}[section]

\newtheorem{lem}[thm]{Lemma}
\newtheorem{prop}[thm]{Proposition}
\newtheorem{defn}[thm]{Definition}

\biboptions{sort&compress,square}
\journal{Journal of Geometry and Physics}
\begin{document}
\begin{frontmatter}
\author{Sining Wei}
\ead{weisn835@nenu.edu.cn}
\author{Yong Wang\corref{cor2}}
\ead{wangy581@nenu.edu.cn}
\cortext[cor2]{Corresponding author.}

\address{School of Mathematics and Statistics, Northeast Normal University,
Changchun, 130024, P.R.China}

\title{Modified Novikov Operators and the Kastler-Kalau-Walze type theorem for manifolds with boundary }
\begin{abstract}
In this paper, we give two Lichnerowicz type formulas for modified Novikov operators. We prove Kastler-Kalau-Walze type theorems for modified Novikov operators on compact manifolds with (resp.without) boundary. We also compute the spectral action for Witten deformation on 4-dimensional compact manifolds.
\end{abstract}
\end{frontmatter}
\section{Introduction}
 As been well known, the noncommutative residue plays a prominent role in noncommutative geometry which is found in \cite{Gu,Wo}. For this reason, it has been studied extensively by geometers. Connes used the noncommutative residue to derive a conformal 4-dimensional Polyakov action analogy in \cite{Co1}. Connes showed us that the noncommutative residue on a compact manifold $M$ coincided with the Dixmier's trace on pseudodifferential operators of order $-{\rm {dim}}M$ in \cite{Co2}.
Connes has also observed that the noncommutative residue of the square of the inverse of the Dirac
operator was proportional to the Einstein-Hilbert action, which is called the Kastler-Kalau-Walze type theorem now. Kastler \cite{Ka} gave a
brute-force proof of this theorem. Kalau and Walze proved this theorem in the normal coordinates system simultaneously in \cite{KW} .
Ackermann proved that
the Wodzicki residue  of the square of the inverse of the Dirac operator ${\rm  Wres}(D^{-2})$ in turn is essentially the second coefficient
of the heat kernel expansion of $D^{2}$ in \cite{Ac}.

On the other hand, Wang generalized the Connes' results to the case of manifolds with boundary in \cite{Wa1,Wa2},
and proved the Kastler-Kalau-Walze type theorem for the Dirac operator and the signature operator on lower-dimensional manifolds
with boundary \cite{Wa3}. In \cite{Wa3,Wa4}, Wang computed $\widetilde{{\rm Wres}}[\pi^+D^{-1}\circ\pi^+D^{-1}]$ and $\widetilde{{\rm Wres}}[\pi^+D^{-2}\circ\pi^+D^{-2}]$, where the two operators are symmetric, in these cases the boundary term vanished. But for $\widetilde{{\rm Wres}}[\pi^+D^{-1}\circ\pi^+D^{-3}]$, Wang got a nonvanishing boundary term \cite{Wa5}, and give a theoretical explanation for gravitational action on boundary. In others words, Wang provides a kind of method to study the Kastler-Kalau-Walze type theorem for manifolds with boundary.
In \cite{lkl}, L\'{o}pez and his collaborators introduced an elliptic differential operator which is called Novikov operator. The motivation of this paper is
to prove the Kastler-Kalau-Walze type theorem for Novikov operators on manifolds with boundary.
In \cite{LI}, Iochum and Levy computed heat kernel coefficients for Dirac operators with one-form perturbations and proved that there are no tadpoles for compact spin manifolds without boundary. In \cite{AA}, Sitarz and Zajac investigated the spectral action for scalar perturbations of Dirac operators. In \cite{SPH}, Hanisch and his collaborators derived a formula for the gravitational part of the spectral action for Dirac operators on 4-dimensional spin manifolds with totally antisymmetric torsion. In \cite{wpz}, Zhang introduced an elliptic differential operator which is called Witten deformation. Motivated by \cite{wpz}, \cite{AA}, \cite{LI} and \cite{SPH}, we will compute the spectral action for Witten deformation on 4-dimensional compact manifolds in this paper.

The framework of this paper is organized as follows. Firstly, in Section 2, we give the definition of modified Novikov operators and the Lichnerowicz formulas associate to modified Novikov operators. We study the symbols of some operators associate to modified Novikov operators, by using symbols of operators associate to modified Novikov operators, we can prove the Kastler-Kalau-Walze type theorem for manifolds with boundary in Section 3 and in Section 4. In Section 5, we compute the spectral action for Witten deformation on 4-dimensional compact manifolds.

%%%%%%%%%%第二部分%%%
\section{Modified Novikov Operators and its Lichnerowicz formula}

In this section, we firstly recall the definition of Novikov Operator (see details in [16]). Let $M$ be a $n$-dimensional ($n\geq 3$) oriented compact Riemannian manifold with a Riemannian metric $g^{M}$. The de Rham derivative $d$ is an elliptic differential operator on $C^\infty(M;\wedge^*T^*M)$. Then we have the de Rham coderivative $\delta=d^*$, the symmetric operators $D=d+\delta$ and $\Delta=D^2=d\delta+\delta d$ (the Laplacian).

With more generality, we take any closed $\theta \in C^\infty(M;T^*M)$. For the sake of simplicity, we assume that $\theta$ is real. Let $V\in \mathfrak{X}(M)$ be determined by $g(V,\cdot)=\theta$. Then we have the Novikov operators defined by $\theta$, depending on $z\in \mathbb{C}$ in [16],
\begin{eqnarray*}
d_z&=&d+z(\theta\wedge),~~\delta_z=d_z^*=\delta+\overline{z}(\theta\wedge)^*,\nonumber\\
D_z&=&d_z+\delta_z=(d+\delta)+z(\theta\wedge)+\overline{z}(\theta\wedge)^*\nonumber\\
   &=&(d+\delta)+[Rez(\theta\wedge)+Rez(\theta\wedge)^*]
                +i[Imz(\theta\wedge)-Imz(\theta\wedge)^*]\nonumber\\
   &=&(d+\delta)+Rez[\theta\wedge+(\theta\wedge)^*]
                +iImz[\theta\wedge-(\theta\wedge)^*]\nonumber\\
   &=&(d+\delta)+Rez\bar{c}(\theta)
                +iImzc(\theta),
\end{eqnarray*}
where $Rez$ is the real part of $z$, $Imz$ is the imaginary part of $z$, $\bar{c}(\theta)=(\theta)^*\wedge+(\theta\wedge)^*$,
$c(\theta)=(\theta)^*\wedge-(\theta\wedge)^*$.

In this paper, we consider the modified Novikov operators, for $\theta,~\theta'\in \Gamma(TM)$ we define that
\begin{eqnarray*}
\widehat{D}&=&d+\delta+\bar{c}(\theta)+c(\theta'),~~\widehat{D}^*=d+\delta+\bar{c}(\theta)-c(\theta').
\end{eqnarray*}
where $\bar{c}(\theta)=(\theta)^*\wedge+(\theta\wedge)^*$,
$c(\theta')=(\theta')^*\wedge-(\theta'\wedge)^*$, where $\theta^*=g(\theta,\cdot)$, $(\theta')^*=g(\theta',\cdot)$.

Let $\nabla^L$ be the Levi-Civita connection about $g^M$. In the local coordinates $\{x_i; 1\leq i\leq n\}$ and the
fixed orthonormal frame $\{\widetilde{e_1},\cdots,\widetilde{e_n}\}$, the connection matrix $(\omega_{s,t})$ is defined by
\begin{equation}
\nabla^L(\widetilde{e_1},\cdots,\widetilde{e_n})= (\widetilde{e_1},\cdots,\widetilde{e_n})(\omega_{s,t}).
\end{equation}
 Let $\epsilon (\widetilde{e_j*} ),~\iota (\widetilde{e_j*} )$ be the exterior and interior multiplications respectively and $c(\widetilde{e_j})$ be the Clifford action.
Suppose that $\partial_{i}$ is a natural local frame on $TM$
and $(g^{ij})_{1\leq i,j\leq n}$ is the inverse matrix associated to the metric
matrix  $(g_{ij})_{1\leq i,j\leq n}$ on $M$. Write
\begin{equation}
c(\widetilde{e_j})=\epsilon (\widetilde{e_j*} )-\iota (\widetilde{e_j*} );~~
\bar{c}(\widetilde{e_j})=\epsilon (\widetilde{e_j*} )+\iota
(\widetilde{e_j*} ).
\end{equation}
 The modified Novikov Operators $\widehat{D}$ and $\widehat{D}^*$ are defined by
\begin{equation}
\widehat{D}=d+\delta+\bar{c}(\theta)+c(\theta')=\sum^n_{i=1}c(\widetilde{e_i})\bigg[\widetilde{e_i}+\frac{1}{4}\sum_{s,t}\omega_{s,t}
(\widetilde{e_i})[\bar{c}(\widetilde{e_s})\bar{c}(\widetilde{e_t})
-c(\widetilde{e_s})c(\widetilde{e_t})]\bigg]+\bar{c}(\theta)+c(\theta');
\end{equation}
\begin{equation}
\widehat{D}^*=d+\delta+\bar{c}(\theta)-c(\theta')=\sum^n_{i=1}c(\widetilde{e_i})\bigg[\widetilde{e_i}+\frac{1}{4}\sum_{s,t}\omega_{s,t}
(\widetilde{e_i})[\bar{c}(\widetilde{e_s})\bar{c}(\widetilde{e_t})
-c(\widetilde{e_s})c(\widetilde{e_t})]\bigg]+\bar{c}(\theta)-c(\theta').
\end{equation}

We first establish the main theorem in this section. One has the following Lichneriowicz formulas.

\begin{thm} The following equalities hold:
\begin{eqnarray}
\widehat{D}^*\widehat{D}&=&-\Big[g^{ij}(\nabla_{\partial_{i}}\nabla_{\partial_{j}}-
\nabla_{\nabla^{L}_{\partial_{i}}\partial_{j}})\Big]
-\frac{1}{8}\sum_{ijkl}R_{ijkl}\bar{c}(\widetilde{e_i})\bar{c}(\widetilde{e_j})
c(\widetilde{e_k})c(\widetilde{e_l})
+\sum_{i}c(\widetilde{e_i})\bar{c}(\nabla_{\widetilde{e_i}}^{TM}\theta)
+\frac{1}{4}s\nonumber\\
&&-c(\theta')\bar{c}(\theta)
+\bar{c}(\theta)c(\theta')+|\theta|^2+|\theta'|^2+\frac{1}{4}\sum_{i}
[c(\widetilde{e_{i}})c(\theta')-c(\theta')c(\widetilde{e_{i}})]^2
-g(\widetilde{e_{j}}
,\nabla^{TM}_{\widetilde{e_{j}}}\theta'),\\
\widehat{D}^2&=&-\Big[g^{ij}(\nabla_{\partial_{i}}\nabla_{\partial_{j}}-
\nabla_{\nabla^{L}_{\partial_{i}}\partial_{j}})\Big]
-\frac{1}{8}\sum_{ijkl}R_{ijkl}\bar{c}(\widetilde{e_i})\bar{c}(\widetilde{e_j})
c(\widetilde{e_k})c(\widetilde{e_l})
+\sum_{i}c(\widetilde{e_i})\bar{c}(\nabla_{\widetilde{e_i}}^{TM}\theta)
\nonumber\\
&&+\frac{1}{4}s+c(\theta')\bar{c}(\theta)
+\bar{c}(\theta)c(\theta')+|\theta|^2-|\theta'|^2+\frac{1}{4}\sum_{i}
[c(\widetilde{e_{i}})c(\theta')+c(\theta')c(\widetilde{e_{i}})]^2\nonumber\\
&&-\frac{1}{2}[c(\nabla^{TM}_{\widetilde{e_{j}}}\theta')c(\widetilde{e_{j}})
-c(\widetilde{e_{j}})c(\nabla^{TM}_{\widetilde{e_{j}}}\theta')],
\end{eqnarray}
where $s$ is the scalar curvature.
\end{thm}

In order to prove Theorem 2.1,
we recall the basic notions of Laplace type operators. Let $M$ be a smooth compact oriented Riemannian $n$-dimensional manifolds without boundary and $V'$ be a vector bundle on $M$. Any differential operator $P$ of Laplace type has locally the form
\begin{equation}
P=-(g^{ij}\partial_i\partial_j+A^i\partial_i+B),
\end{equation}
where $\partial_{i}$  is a natural local frame on $TM$
and $(g^{ij})_{1\leq i,j\leq n}$ is the inverse matrix associated to the metric
matrix  $(g_{ij})_{1\leq i,j\leq n}$ on $M$,
 and $A^{i}$ and $B$ are smooth sections
of $\textrm{End}(V')$ on $M$ (endomorphism). If $P$ is a Laplace type
operator with the form (2.7), then there is a unique
connection $\nabla$ on $V'$ and a unique endomorphism $E$ such that
 \begin{equation}
P=-\Big[g^{ij}(\nabla_{\partial_{i}}\nabla_{\partial_{j}}-
 \nabla_{\nabla^{L}_{\partial_{i}}\partial_{j}})+E\Big],
\end{equation}
where $\nabla^{L}$ is the Levi-Civita connection on $M$. Moreover
(with local frames of $T^{*}M$ and $V'$), $\nabla_{\partial_{i}}=\partial_{i}+\omega_{i} $
and $E$ are related to $g^{ij}$, $A^{i}$ and $B$ through
 \begin{eqnarray}
&&\omega_{i}=\frac{1}{2}g_{ij}\big(A^{i}+g^{kl}\Gamma_{ kl}^{j} \texttt{id}\big),\\
&&E=B-g^{ij}\big(\partial_{i}(\omega_{j})+\omega_{i}\omega_{j}-\omega_{k}\Gamma_{ ij}^{k} \big),
\end{eqnarray}
where $\Gamma_{ kl}^{j}$ is the  Christoffel coefficient of $\nabla^{L}$.

By Proposition 4.6 of \cite{wpz}, we have
\begin{equation}
(d+\delta+\bar{c}(\theta))^2=(d+\delta)^{2}
+\sum_{i}c(\widetilde{e_i})\bar{c}(\nabla_{\widetilde{e_i}}^{TM}\theta)+|\theta|^{2}.
\end{equation}
By \cite{Y}, the local expression of $(d+\delta)^{2}$ is
\begin{equation}
(d+\delta)^{2}
=-\triangle_{0}-\frac{1}{8}\sum_{ijkl}R_{ijkl}\bar{c}(\widetilde{e_i})\bar{c}(\widetilde{e_j})c(\widetilde{e_k})c(\widetilde{e_l})+\frac{1}{4}s.
\end{equation}
Let $g^{ij}=g(dx_{i},dx_{j})$, $\xi=\sum_{k}\xi_{j}dx_{j}$ and $\nabla^L_{\partial_{i}}\partial_{j}=\sum_{k}\Gamma_{ij}^{k}\partial_{k}$,  we denote that
\begin{eqnarray}
&&\sigma_{i}=-\frac{1}{4}\sum_{s,t}\omega_{s,t}
(\widetilde{e_i})c(\widetilde{e_s})c(\widetilde{e_t})
;~~~a_{i}=\frac{1}{4}\sum_{s,t}\omega_{s,t}
(\widetilde{e_i})\bar{c}(\widetilde{e_s})\bar{c}(\widetilde{e_t});\nonumber\\
&&\xi^{j}=g^{ij}\xi_{i};~~~~\Gamma^{k}=g^{ij}\Gamma_{ij}^{k};~~~~\sigma^{j}=g^{ij}\sigma_{i};
~~~~a^{j}=g^{ij}a_{i}.
\end{eqnarray}
Then the modified Novikov Operators $\widehat{D}$ and $\widehat{D}^*$ can be written as
\begin{equation}
\widehat{D}=\sum^n_{i=1}c(\widetilde{e_i})[\widetilde{e_i}+a_{i}+\sigma_{i}]
+\bar{c}(\theta)+c(\theta');
\end{equation}
\begin{equation}
\widehat{D}^*=\sum^n_{i=1}c(\widetilde{e_i})[\widetilde{e_i}+a_{i}+\sigma_{i}]
+\bar{c}(\theta)-c(\theta').
\end{equation}
By \cite{Y} and \cite{Ac}, we have
\begin{equation}
-\triangle_{0}=\Delta=-g^{ij}(\nabla^L_{i}\nabla^L_{j}-\Gamma_{ij}^{k}\nabla^L_{k}).
\end{equation}
We note that
\begin{eqnarray}
\widehat{D}^*\widehat{D}&=&(d+\delta+\bar{c}(\theta))^2+(d+\delta)c(\theta')
+\bar{c}(\theta)c(\theta')-c(\theta')(d+\delta)-c(\theta')\bar{c}(\theta)
+|\theta'|^2,
\end{eqnarray}
\begin{eqnarray}
-c(\theta')(d+\delta)+(d+\delta)c(\theta')&=&\sum_{i,j}g^{i,j}\Big[c(\partial_{i})c(\theta')
-c(\theta')c(\partial_{i})\Big]\partial_{j}
-\sum_{i,j}g^{i,j}\Big[c(\theta')c(\partial_{i})\sigma_{i}+c(\theta')c(\partial_{i})a_{i}\nonumber\\
&&+c(\partial_{i})\partial_{i}(c(\theta'))+c(\partial_{i})\sigma_{i}c(\theta')
+c(\partial_{i})a_{i}c(\theta')\Big],
\end{eqnarray}
then we obtain
\begin{eqnarray}
\widehat{D}^*\widehat{D}&=&-\sum_{i,j}g^{i,j}\Big[\partial_{i}\partial_{j}
+2\sigma_{i}\partial_{j}+2a_{i}\partial_{j}-\Gamma_{i,j}^{k}\partial_{k}
+(\partial_{i}\sigma_{j})
+(\partial_{i}a_{j})
+\sigma_{i}\sigma_{j}+\sigma_{i}a_{j}+a_{i}\sigma_{j}+a_{i}a_{j} -\Gamma_{i,j}^{k}\sigma_{k}\nonumber\\
&&-\Gamma_{i,j}^{k}a_{k}\Big]
+\sum_{i,j}g^{i,j}\Big[c(\partial_{i})c(\theta')-c(\theta')c(\partial_{i})\Big]\partial_{j}
-\sum_{i,j}g^{i,j}\Big[c(\theta')c(\partial_{i})\sigma_{i}+c(\theta')c(\partial_{i})a_{i}
-c(\partial_{i})\partial_{i}(c(\theta'))\nonumber\\
&&-c(\partial_{i})\sigma_{i}c(\theta')-c(\partial_{i})a_{i}c(\theta')\Big]
-\frac{1}{8}\sum_{ijkl}R_{ijkl}\bar{c}(\widetilde{e_i})\bar{c}(\widetilde{e_j})
c(\widetilde{e_k})c(\widetilde{e_l})+\frac{1}{4}s
+\sum_{i}c(\widetilde{e_i})\bar{c}(\nabla_{\widetilde{e_i}}^{TM}\theta)+|\theta|^{2}\nonumber\\
&&+|\theta'|^2+\bar{c}(\theta)c(\theta')-c(\theta')\bar{c}(\theta).
\end{eqnarray}
Similarly, we have
\begin{eqnarray}
\widehat{D}^2&=&-\sum_{i,j}g^{i,j}\Big[\partial_{i}\partial_{j}
+2\sigma_{i}\partial_{j}+2a_{i}\partial_{j}-\Gamma_{i,j}^{k}\partial_{k}
+(\partial_{i}\sigma_{j})
+(\partial_{i}a_{j})
+\sigma_{i}\sigma_{j}+\sigma_{i}a_{j}+a_{i}\sigma_{j}+a_{i}a_{j} -\Gamma_{i,j}^{k}\sigma_{k}\nonumber\\
&&-\Gamma_{i,j}^{k}a_{k}\Big]
+\sum_{i,j}g^{i,j}\Big[c(\partial_{i})c(\theta')+c(\theta')c(\partial_{i})\Big]\partial_{j}
+\sum_{i,j}g^{i,j}\Big[c(\theta')c(\partial_{i})\sigma_{i}+c(\theta')c(\partial_{i})a_{i}
+c(\partial_{i})\partial_{i}(c(\theta'))\nonumber\\
&&+c(\partial_{i})\sigma_{i}c(\theta')+c(\partial_{i})a_{i}c(\theta')\Big]
-\frac{1}{8}\sum_{ijkl}R_{ijkl}\bar{c}(\widetilde{e_i})\bar{c}(\widetilde{e_j})
c(\widetilde{e_k})c(\widetilde{e_l})+\frac{1}{4}s
+\sum_{i}c(\widetilde{e_i})\bar{c}(\nabla_{\widetilde{e_i}}^{TM}\theta)+|\theta|^{2}\nonumber\\
&&-|\theta'|^2+\bar{c}(\theta)c(\theta')+c(\theta')\bar{c}(\theta).
\end{eqnarray}
By (2.8), (2.9), (2.10) and (2.19), we have
\begin{eqnarray}
(\omega_{i})_{\widehat{D}^*\widehat{D}}&=&\sigma_{i}+a_{i}
-\frac{1}{2}\Big[c(\partial_{i})c(\theta')-c(\theta')c(\partial_{i})\Big],\\
E_{\widehat{D}^*\widehat{D}}&=&-c(\partial_{i})\sigma^{i}c(\theta')-c(\partial_{i})a^{i}c(\theta')
+\frac{1}{8}\sum_{ijkl}R_{ijkl}\bar{c}(\widetilde{e_i})\bar{c}(\widetilde{e_j})
c(\widetilde{e_k})c(\widetilde{e_l})
-\sum_{i}c(\widetilde{e_i})\bar{c}(\nabla_{\widetilde{e_i}}^{TM}\theta)
-|\theta|^{2}-|\theta'|^2\nonumber\\
&&-\frac{1}{4}s+c(\theta')\bar{c}(\theta)
+c(\theta')c(\partial_{i})\sigma^{i}+c(\theta')c(\partial_{i})a^{i}
-c(\partial_{i})\partial^{i}(c(\theta'))
+\frac{1}{2}\partial^{j}[c(\partial_{j})c(\theta')
-c(\theta')c(\partial_{j})]\nonumber\\
&&-\frac{1}{2}[c(\partial_{j})c(\theta')
-c(\theta')c(\partial_{j})](\sigma^{j}+a^{j})
-\frac{g^{ij}}{4}[c(\partial_{i})c(\theta')
-c(\theta')c(\partial_{i})][c(\partial_{j})c(\theta')
-c(\theta')c(\partial_{j})]\nonumber\\
&&-\frac{1}{2}\Gamma^{k}[c(\partial_{k})c(\theta')
-c(\theta')c(\partial_{k})]
-\bar{c}(\theta)c(\theta')
-\frac{1}{2}(\sigma^{j}+a^{j})
[c(\partial_{j})c(\theta')
-c(\theta')c(\partial_{j})].
\end{eqnarray}
For a smooth vector field $Y$ on $M$, let $c(Y)$
denote the Clifford action. Since $E$ is globally
defined on $M$, taking normal coordinates at $x_0$, we have
$\sigma^{i}(x_0)=0$, $a^{i}(x_0)=0$, $\partial^{j}[c(\partial_{j})](x_0)=0$,
$\Gamma^k(x_0)=0$, $g^{ij}(x_0)=\delta^j_i$, so that
\begin{eqnarray}
E_{\widehat{D}^*\widehat{D}}(x_0)&=&\frac{1}{8}\sum_{ijkl}R_{ijkl}\bar{c}(\widetilde{e_i})\bar{c}(\widetilde{e_j})
c(\widetilde{e_k})c(\widetilde{e_l})
-\sum_{i}c(\widetilde{e_i})\bar{c}(\nabla_{\widetilde{e_i}}^{TM}\theta)
-\frac{1}{4}s+c(\theta')\bar{c}(\theta)-\bar{c}(\theta)c(\theta')
-|\theta|^2-|\theta'|^2\nonumber\\
&&-\frac{1}{4}\sum_{i}[c(\widetilde{e_{i}})c(\theta')
-c(\theta')c(\widetilde{e_{i}})]^2
-\frac{1}{2}[c(\widetilde{e_{j}})\widetilde{e_{j}}(c(\theta'))
+\widetilde{e_{j}}(c(\theta'))c(\widetilde{e_{j}})]\nonumber\\
&=&\frac{1}{8}\sum_{ijkl}R_{ijkl}\bar{c}(\widetilde{e_i})\bar{c}(\widetilde{e_j})
c(\widetilde{e_k})c(\widetilde{e_l})
-\sum_{i}c(\widetilde{e_i})\bar{c}(\nabla_{\widetilde{e_i}}^{TM}\theta)
-\frac{1}{4}s+c(\theta')\bar{c}(\theta)-\bar{c}(\theta)c(\theta')
-|\theta|^2-|\theta'|^2\nonumber\\
&&-\frac{1}{4}\sum_{i}[c(\widetilde{e_{i}})c(\theta')-c(\theta')c(\widetilde{e_{i}})]^2
+g(\widetilde{e_{j}}
,\nabla^{TM}_{\widetilde{e_{j}}}\theta').
\end{eqnarray}
Similarly, we have
\begin{eqnarray}
E_{\widehat{D}^2}(x_0)
&=&\frac{1}{8}\sum_{ijkl}R_{ijkl}\bar{c}(\widetilde{e_i})\bar{c}(\widetilde{e_j})
c(\widetilde{e_k})c(\widetilde{e_l})
-\sum_{i}c(\widetilde{e_i})\bar{c}(\nabla_{\widetilde{e_i}}^{TM}\theta)
-\frac{1}{4}s-c(\theta')\bar{c}(\theta)-\bar{c}(\theta)c(\theta')
-|\theta|^2+|\theta'|^2\nonumber\\
&&-\frac{1}{4}\sum_{i}[c(\widetilde{e_{i}})c(\theta')+c(\theta')c(\widetilde{e_{i}})]^2
+\frac{1}{2}[c(\nabla^{TM}_{\widetilde{e_{j}}}\theta')c(\widetilde{e_{j}})
-c(\widetilde{e_{j}})c(\nabla^{TM}_{\widetilde{e_{j}}}\theta')],
\end{eqnarray}
which, together with (2.7), yields Theorem 2.1.

The non-commutative residue of a generalized laplacian $\widetilde{\Delta}$ is expressed as by \cite{Ac}

\begin{equation}
(n-2)\Phi_{2}(\widetilde{\Delta})=(4\pi)^{-\frac{n}{2}}\Gamma(\frac{n}{2})\widetilde{res}(\widetilde{\Delta}^{-\frac{n}{2}+1}),
\end{equation}
where $\Phi_{2}(\widetilde{\Delta})$ denotes the integral over the diagonal part of the second
coefficient of the heat kernel expansion of $\widetilde{\Delta}$.
Now let $\widetilde{\Delta}=\widehat{D}^*\widehat{D}$. Since $\widehat{D}^*\widehat{D}$ is a generalized laplacian , we can suppose $\widehat{D}^*\widehat{D}=\Delta-E$, then, we have
\begin{eqnarray}
{\rm Wres}(\widehat{D}^*\widehat{D})^{-\frac{n-2}{2}}
=\frac{(n-2)(4\pi)^{\frac{n}{2}}}{(\frac{n}{2}-1)!}\int_{M}{\rm tr}(\frac{1}{6}s+E_{\widehat{D}^*\widehat{D}})d{\rm Vol_{M} },
\end{eqnarray}
where ${\rm Wres}$ denote the noncommutative residue.

Similarly, we have
\begin{eqnarray}
{\rm Wres}(\widehat{D}^2)^{-\frac{n-2}{2}}
=\frac{(n-2)(4\pi)^{\frac{n}{2}}}{(\frac{n}{2}-1)!}\int_{M}{\rm tr}(\frac{1}{6}s+E_{\widehat{D}^2})d{\rm Vol_{M} },
\end{eqnarray}
where ${\rm Wres}$ denote the noncommutative residue.

\begin{thm} For even $n$-dimensional compact oriented manifolds without boundary, the following equalities holds:
\begin{eqnarray}
{\rm Wres}(\widehat{D}^*\widehat{D})^{-\frac{n-2}{2}}
&=&\frac{(n-2)(4\pi)^{\frac{n}{2}}}{(\frac{n}{2}-1)!}\int_{M}
2^n\bigg(-\frac{1}{12}s-|\theta|^2+(n-2)|\theta'|^2+    g(\widetilde{e_{j}},\nabla^{TM}_{\widetilde{e_{j}}}\theta')\bigg)d{\rm Vol_{M}},\\
{\rm Wres}(\widehat{D}^2)^{-\frac{n-2}{2}}
&=&\frac{(n-2)(4\pi)^{\frac{n}{2}}}{(\frac{n}{2}-1)!}\int_{M}2^n\bigg(
-\frac{1}{12}s-|\theta|^2\bigg)d{\rm Vol_{M}}.
\end{eqnarray}
where $s$ is the scalar curvature.
\end{thm}

%%%%%%%%%第3部分%%%%%%%%
\section{A Kastler-Kalau-Walze type theorem for $4$-dimensional manifolds with boundary}
In this section, we prove the Kastler-Kalau-Walze type theorem for $4$-dimensional oriented compact manifold with boundary. We firstly give some basic facts and formulas about Boutet de
Monvel's calculus and the definition of the noncommutative residue for manifolds with boundary (see details in Section 2 in \cite{Wa3}).

Let $U\subset M$ be a collar neighborhood of $\partial M$ which is diffeomorphic with $\partial M\times [0,1)$. By the definition of $h(x_n)\in C^{\infty}([0,1))$
and $h(x_n)>0$, there exists $\tilde{h}\in C^{\infty}((-\varepsilon,1))$ such that $\tilde{h}|_{[0,1)}=h$ and $\tilde{h}>0$ for some
sufficiently small $\varepsilon>0$. Then there exists a metric $\widehat{g}$ on $\widehat{M}=M\bigcup_{\partial M}\partial M\times
(-\varepsilon,0]$ which has the form on $U\bigcup_{\partial M}\partial M\times (-\varepsilon,0 ]$
\begin{equation}
\widehat{g}=\frac{1}{\tilde{h}(x_{n})}g^{\partial M}+dx _{n}^{2} ,
\end{equation}
such that $\widehat{g}|_{M}=g$. We fix a metric $\widehat{g}$ on the $\widehat{M}$ such that $\widehat{g}|_{M}=g$.

Let  \begin{equation*}
F:L^2({\bf R}_t)\rightarrow L^2({\bf R}_v);~F(u)(v)=\int e^{-ivt}u(t)dt
\end{equation*}
denote the Fourier transformation and
$\Phi(\overline{{\bf R}^+}) =r^+\Phi({\bf R})$ (similarly define $\Phi(\overline{{\bf R}^-}$)), where $\Phi({\bf R})$
denotes the Schwartz space and
\begin{equation*}
r^{+}:C^\infty ({\bf R})\rightarrow C^\infty (\overline{{\bf R}^+});~ f\rightarrow f|\overline{{\bf R}^+};~
\overline{{\bf R}^+}=\{x\geq0;x\in {\bf R}\}.
\end{equation*}
We define $H^+=F(\Phi(\overline{{\bf R}^+}));~ H^-_0=F(\Phi(\overline{{\bf R}^-}))$ which are orthogonal to each other. We have the following
 property: $h\in H^+~(H^-_0)$ if and only if $h\in C^\infty({\bf R})$ which has an analytic extension to the lower (upper) complex
half-plane $\{{\rm Im}\xi<0\}~(\{{\rm Im}\xi>0\})$ such that for all nonnegative integer $l$,
 \begin{equation*}
\frac{d^{l}h}{d\xi^l}(\xi)\sim\sum^{\infty}_{k=1}\frac{d^l}{d\xi^l}(\frac{c_k}{\xi^k}),
\end{equation*}
as $|\xi|\rightarrow +\infty,{\rm Im}\xi\leq0~({\rm Im}\xi\geq0)$.

 Let $H'$ be the space of all polynomials and $H^-=H^-_0\bigoplus H';~H=H^+\bigoplus H^-.$ Denote by $\pi^+~(\pi^-)$ respectively the
 projection on $H^+~(H^-)$. For calculations, we take $H=\widetilde H=\{$rational functions having no poles on the real axis$\}$ ($\tilde{H}$
 is a dense set in the topology of $H$). Then on $\tilde{H}$,
 \begin{equation}
\pi^+h(\xi_0)=\frac{1}{2\pi i}\lim_{u\rightarrow 0^{-}}\int_{\Gamma^+}\frac{h(\xi)}{\xi_0+iu-\xi}d\xi,
\end{equation}
where $\Gamma^+$ is a Jordan close curve
included ${\rm Im}(\xi)>0$ surrounding all the singularities of $h$ in the upper half-plane and
$\xi_0\in {\bf R}$. Similarly, define $\pi'$ on $\tilde{H}$,
\begin{equation}
\pi'h=\frac{1}{2\pi}\int_{\Gamma^+}h(\xi)d\xi.
\end{equation}
So, $\pi'(H^-)=0$. For $h\in H\bigcap L^1({\bf R})$, $\pi'h=\frac{1}{2\pi}\int_{{\bf R}}h(v)dv$ and for $h\in H^+\bigcap L^1({\bf R})$, $\pi'h=0$.

Let $M$ be a $n$-dimensional compact oriented manifold with boundary $\partial M$.
Denote by $\mathcal{B}$ Boutet de Monvel's algebra, we recall the main theorem in \cite{FGLS,Wa3}.
\begin{thm}\label{th:32}{\rm\cite{FGLS}}{\bf(Fedosov-Golse-Leichtnam-Schrohe)}
 Let $X$ and $\partial X$ be connected, ${\rm dim}X=n\geq3$,
 $A=\left(\begin{array}{lcr}\pi^+P+G &   K \\
T &  S    \end{array}\right)$ $\in \mathcal{B}$ , and denote by $p$, $b$ and $s$ the local symbols of $P,G$ and $S$ respectively.
 Define:
 \begin{eqnarray}
{\rm{\widetilde{Wres}}}(A)&=&\int_X\int_{\bf S}{\rm{tr}}_E\left[p_{-n}(x,\xi)\right]\sigma(\xi)dx \nonumber\\
&&+2\pi\int_ {\partial X}\int_{\bf S'}\left\{{\rm tr}_E\left[({\rm{tr}}b_{-n})(x',\xi')\right]+{\rm{tr}}
_F\left[s_{1-n}(x',\xi')\right]\right\}\sigma(\xi')dx',
\end{eqnarray}
Then~~ a) ${\rm \widetilde{Wres}}([A,B])=0 $, for any
$A,B\in\mathcal{B}$;~~ b) It is a unique continuous trace on
$\mathcal{B}/\mathcal{B}^{-\infty}$.
\end{thm}

\begin{defn}{\rm\cite{Wa3} }
Lower dimensional volumes of spin manifolds with boundary are defined by
 \begin{equation}
{\rm Vol}^{(p_1,p_2)}_nM:= \widetilde{{\rm Wres}}[\pi^+D^{-p_1}\circ\pi^+D^{-p_2}],
\end{equation}
\end{defn}
 By (2.1.4)-(2.1.8) in \cite{Wa3}, we get
\begin{equation}
\widetilde{{\rm Wres}}[\pi^+D^{-p_1}\circ\pi^+D^{-p_2}]=\int_M\int_{|\xi|=1}{\rm
trace}_{\wedge^*T^*M}[\sigma_{-n}(D^{-p_1-p_2})]\sigma(\xi)dx+\int_{\partial M}\Phi,
\end{equation}
and
\begin{eqnarray}
\Phi &=&\int_{|\xi'|=1}\int^{+\infty}_{-\infty}\sum^{\infty}_{j, k=0}\sum\frac{(-i)^{|\alpha|+j+k+1}}{\alpha!(j+k+1)!}
\times {\rm trace}_{\wedge^*T^*M}[\partial^j_{x_n}\partial^\alpha_{\xi'}\partial^k_{\xi_n}\sigma^+_{r}(D^{-p_1})(x',0,\xi',\xi_n)
\nonumber\\
&&\times\partial^\alpha_{x'}\partial^{j+1}_{\xi_n}\partial^k_{x_n}\sigma_{l}(D^{-p_2})(x',0,\xi',\xi_n)]d\xi_n\sigma(\xi')dx',
\end{eqnarray}
 where the sum is taken over $r+l-k-|\alpha|-j-1=-n,~~r\leq -p_1,l\leq -p_2$.

 Since $[\sigma_{-n}(D^{-p_1-p_2})]|_M$ has the same expression as $\sigma_{-n}(D^{-p_1-p_2})$ in the case of manifolds without
boundary, so locally we can use the computations \cite{Ka}, \cite{KW}, \cite{Wa3}, \cite{Po} to compute the first term.
%The following proposition is the motivation
%of the definition of lower dimensional volumes of spin manifolds with boundary \cite{Wa4}.

For any fixed point $x_0\in\partial M$, we choose the normal coordinates
$U$ of $x_0$ in $\partial M$ (not in $M$) and compute $\Phi(x_0)$ in the coordinates $\widetilde{U}=U\times [0,1)\subset M$ and the
metric $\frac{1}{h(x_n)}g^{\partial M}+dx_n^2.$ The dual metric of $g^M$ on $\widetilde{U}$ is ${h(x_n)}g^{\partial M}+dx_n^2.$  Write
$g^M_{ij}=g^M(\frac{\partial}{\partial x_i},\frac{\partial}{\partial x_j});~ g_M^{ij}=g^M(dx_i,dx_j)$, then

\begin{equation}
[g^M_{i,j}]= \left[\begin{array}{lcr}
  \frac{1}{h(x_n)}[g_{i,j}^{\partial M}]  & 0  \\
   0  &  1
\end{array}\right];~~~
[g_M^{i,j}]= \left[\begin{array}{lcr}
  h(x_n)[g^{i,j}_{\partial M}]  & 0  \\
   0  &  1
\end{array}\right],
\end{equation}
and
\begin{equation}
\partial_{x_s}g_{ij}^{\partial M}(x_0)=0, 1\leq i,j\leq n-1; ~~~g_{ij}^M(x_0)=\delta_{ij}.
\end{equation}
We will give following three lemmas as computation tools.
\begin{lem}{\rm \cite{Wa3}}\label{le:32}
With the metric $g^{M}$ on $M$ near the boundary
\begin{eqnarray}
\partial_{x_j}(|\xi|_{g^M}^2)(x_0)&=&\left\{
       \begin{array}{c}
        0,  ~~~~~~~~~~ ~~~~~~~~~~ ~~~~~~~~~~~~~{\rm if }~j<n, \\[2pt]
       h'(0)|\xi'|^{2}_{g^{\partial M}},~~~~~~~~~~~~~~~~~~~~{\rm if }~j=n;
       \end{array}
    \right. \\
\partial_{x_j}[c(\xi)](x_0)&=&\left\{
       \begin{array}{c}
      0,  ~~~~~~~~~~ ~~~~~~~~~~ ~~~~~~~~~~~~~{\rm if }~j<n,\\[2pt]
\partial x_{n}(c(\xi'))(x_{0}), ~~~~~~~~~~~~~~~~~{\rm if }~j=n,
       \end{array}
    \right.
\end{eqnarray}
where $\xi=\xi'+\xi_{n}dx_{n}$.
\end{lem}
\begin{lem}{\rm \cite{Wa3}}\label{le:32}With the metric $g^{M}$ on $M$ near the boundary
\begin{eqnarray}
\omega_{s,t}(\widetilde{e_i})(x_0)&=&\left\{
       \begin{array}{c}
        \omega_{n,i}(\widetilde{e_i})(x_0)=\frac{1}{2}h'(0),  ~~~~~~~~~~ ~~~~~~~~~~~{\rm if }~s=n,t=i,i<n; \\[2pt]
       \omega_{i,n}(\widetilde{e_i})(x_0)=-\frac{1}{2}h'(0),~~~~~~~~~~~~~~~~~~~{\rm if }~s=i,t=n,i<n;\\[2pt]
    \omega_{s,t}(\widetilde{e_i})(x_0)=0,~~~~~~~~~~~~~~~~~~~~~~~~~~~other~cases~~~~~~~~~,
       \end{array}
    \right.
\end{eqnarray}
where $(\omega_{s,t})$ denotes the connection matrix of Levi-Civita connection $\nabla^L$.
\end{lem}
\begin{lem}{\rm \cite{Wa3}}~{\it When $i<n$, then }
$$\Gamma^n_{ii}(x_0)=\frac{1}{2}h'(0);~\Gamma^i_{ni}(x_0)=-\frac{1}{2}h'(0);~\Gamma^i_{in}(x_0)=-\frac{1}{2}h'(0),$$
\noindent {\it in other cases}, $ \Gamma_{st}^i(x_0)=0.$
\end{lem}

By (3.6) and (3.7), we firstly compute
\begin{equation}
\widetilde{{\rm Wres}}[\pi^+\widehat{D}^{-1}\circ\pi^+(\widehat{D}^*)^{-1}]=\int_M\int_{|\xi|=1}{\rm
trace}_{\wedge^*T^*M}[\sigma_{-4}((\widehat{D}^*\widehat{D})^{-1})]\sigma(\xi)dx+\int_{\partial M}\Phi,
\end{equation}
where
\begin{eqnarray}
\Phi &=&\int_{|\xi'|=1}\int^{+\infty}_{-\infty}\sum^{\infty}_{j, k=0}\sum\frac{(-i)^{|\alpha|+j+k+1}}{\alpha!(j+k+1)!}
\times {\rm trace}_{\wedge^*T^*M}[\partial^j_{x_n}\partial^\alpha_{\xi'}\partial^k_{\xi_n}\sigma^+_{r}(\widehat{D}^{-1})(x',0,\xi',\xi_n)
\nonumber\\
&&\times\partial^\alpha_{x'}\partial^{j+1}_{\xi_n}\partial^k_{x_n}\sigma_{l}((\widehat{D}^*)^{-1})(x',0,\xi',\xi_n)]d\xi_n\sigma(\xi')dx',
\end{eqnarray}
and the sum is taken over $r+l-k-j-|\alpha|=-3,~~r\leq -1,l\leq-1$.\\

Locally we can use Theorem 2.2 (2.28) to compute the interior of $\widetilde{{\rm Wres}}[\pi^+\widehat{D}^{-1}\circ\pi^+(\widehat{D}^*)^{-1}]$, we have
\begin{eqnarray}
&&\int_M\int_{|\xi|=1}{\rm
trace}_{\wedge^*T^*M}[\sigma_{-4}((\widehat{D}^*\widehat{D})^{-1})]\sigma(\xi)dx\nonumber\\
&=&32\pi^2\int_{M}
\bigg[16g(\widetilde{e_{j}},\nabla^{TM}_{\widetilde{e_{j}}}\theta')-\frac{4}{3}s
-16|\theta|^2+32|\theta'|^2\bigg]d{\rm Vol_{M}}.
\end{eqnarray}

So we only need to compute $\int_{\partial M} \Phi$. Let us now turn to compute the symbols of some operators.
By (2.13)-(2.20), some operators have the following symbols.
\begin{lem} The following identities hold:
\begin{eqnarray}
\sigma_1(\widehat{D})&=&\sigma_1(\widehat{D}^*)=ic(\xi); \nonumber\\ \sigma_0(\widehat{D})&=&\frac{1}{4}\sum_{i,s,t}\omega_{s,t}(\widetilde{e_i})c(\widetilde{e_i})\bar{c}(\widetilde{e_s})\bar{c}(\widetilde{e_t})
-\frac{1}{4}\sum_{i,s,t}\omega_{s,t}(\widetilde{e_i})c(\widetilde{e_i})
c(\widetilde{e_s})c(\widetilde{e_t})+\bar{c}(\theta)+c(\theta'); \nonumber\\
\sigma_0(\widehat{D}^*)&=&\frac{1}{4}\sum_{i,s,t}\omega_{s,t}(\widetilde{e_i})c(\widetilde{e_i})\bar{c}(\widetilde{e_s})\bar{c}(\widetilde{e_t})
-\frac{1}{4}\sum_{i,s,t}\omega_{s,t}(\widetilde{e_i})c(\widetilde{e_i})
c(\widetilde{e_s})c(\widetilde{e_t})+\bar{c}(\theta)-c(\theta').
\end{eqnarray}
\end{lem}

Write
 \begin{eqnarray}
D_x^{\alpha}&=&(-i)^{|\alpha|}\partial_x^{\alpha};
~\sigma(\widehat{D})=p_1+p_0;
~\sigma(\widehat{D}^{-1})=\sum^{\infty}_{j=1}q_{-j}.
\end{eqnarray}

By the composition formula of pseudodifferential operators, we have
\begin{eqnarray}
1=\sigma(\widehat{D}\circ \widehat{D}^{-1})&=&\sum_{\alpha}\frac{1}{\alpha!}\partial^{\alpha}_{\xi}[\sigma(\widehat{D})]
D^{\alpha}_{x}[\sigma(\widehat{D}^{-1})]\nonumber\\
&=&(p_1+p_0)(q_{-1}+q_{-2}+q_{-3}+\cdots)\nonumber\\
& &~~~+\sum_j(\partial_{\xi_j}p_1+\partial_{\xi_j}p_0)(
D_{x_j}q_{-1}+D_{x_j}q_{-2}+D_{x_j}q_{-3}+\cdots)\nonumber\\
&=&p_1q_{-1}+(p_1q_{-2}+p_0q_{-1}+\sum_j\partial_{\xi_j}p_1D_{x_j}q_{-1})+\cdots,
\end{eqnarray}
so
\begin{equation}
q_{-1}=p_1^{-1};~q_{-2}=-p_1^{-1}[p_0p_1^{-1}+\sum_j\partial_{\xi_j}p_1D_{x_j}(p_1^{-1})].
\end{equation}

By Lemma 3.6, we have some symbols of operators.
\begin{lem} The following identities hold:
\begin{eqnarray}
\sigma_{-1}(\widehat{D}^{-1})&=&\sigma_{-1}((\widehat{D}^*)^{-1})=\frac{ic(\xi)}{|\xi|^2};\nonumber\\
\sigma_{-2}(\widehat{D}^{-1})&=&\frac{c(\xi)\sigma_{0}(\widehat{D}^{-1})c(\xi)}{|\xi|^4}+\frac{c(\xi)}{|\xi|^6}\sum_jc(dx_j)
\Big[\partial_{x_j}(c(\xi))|\xi|^2-c(\xi)\partial_{x_j}(|\xi|^2)\Big] ;\nonumber\\
\sigma_{-2}((\widehat{D}^*)^{-1})&=&\frac{c(\xi)\sigma_{0}((\widehat{D}^*)^{-1})c(\xi)}{|\xi|^4}+\frac{c(\xi)}{|\xi|^6}\sum_jc(dx_j)
\Big[\partial_{x_j}(c(\xi))|\xi|^2-c(\xi)\partial_{x_j}(|\xi|^2)\Big]. \end{eqnarray}
\end{lem}

From the remark above, now we can compute $\Phi$ (see formula (3.14) for the definition of $\Phi$). We use ${\rm tr}$ as shorthand of ${\rm trace}$. Since $n=4$, then ${\rm tr}_{\wedge^*T^*M}[{\rm \texttt{id}}]={\rm dim}(\wedge^*(4))=16$, since the sum is taken over $
r+l-k-j-|\alpha|=-3,~~r\leq -1,l\leq-1,$ then we have the following five cases:

~\\
\noindent  {\bf case a)~I)}~$r=-1,~l=-1,~k=j=0,~|\alpha|=1$\\

\noindent By (3.14), we get
\begin{equation}
{\rm case~a)~I)}=-\int_{|\xi'|=1}\int^{+\infty}_{-\infty}\sum_{|\alpha|=1}
 {\rm tr}[\partial^\alpha_{\xi'}\pi^+_{\xi_n}\sigma_{-1}(\widehat{D}^{-1})\times
 \partial^\alpha_{x'}\partial_{\xi_n}\sigma_{-1}((\widehat{D}^*)^{-1})](x_0)d\xi_n\sigma(\xi')dx'.
\end{equation}
By Lemma 3.3, for $i<n$, then
\begin{equation}\partial_{x_i}\left(\frac{ic(\xi)}{|\xi|^2}\right)(x_0)=
\frac{i\partial_{x_i}[c(\xi)](x_0)}{|\xi|^2}
-\frac{ic(\xi)\partial_{x_i}(|\xi|^2)(x_0)}{|\xi|^4}=0,
\end{equation}
\noindent so {\rm case~a)~I)} vanishes.\\

 \noindent  {\bf case a)~II)}~$r=-1,~l=-1,~k=|\alpha|=0,~j=1$\\

\noindent By (3.14), we get
\begin{equation}
{\rm case~a)~II)}=-\frac{1}{2}\int_{|\xi'|=1}\int^{+\infty}_{-\infty} {\rm
trace} [\partial_{x_n}\pi^+_{\xi_n}\sigma_{-1}(\widehat{D}^{-1})\times
\partial_{\xi_n}^2\sigma_{-1}((\widehat{D}^*)^{-1})](x_0)d\xi_n\sigma(\xi')dx'.
\end{equation}
\noindent By Lemma 3.7, we have\\
\begin{eqnarray}\partial^2_{\xi_n}\sigma_{-1}((\widehat{D}^*)^{-1})(x_0)=i\left(-\frac{6\xi_nc(dx_n)+2c(\xi')}
{|\xi|^4}+\frac{8\xi_n^2c(\xi)}{|\xi|^6}\right);
\end{eqnarray}
\begin{eqnarray}
\partial_{x_n}\sigma_{-1}(\widehat{D}^{-1})(x_0)=\frac{i\partial_{x_n}c(\xi')(x_0)}{|\xi|^2}-\frac{ic(\xi)|\xi'|^2h'(0)}{|\xi|^4}.
\end{eqnarray}
By (3.2), (3.3) and the Cauchy integral formula we have
\begin{eqnarray}
\pi^+_{\xi_n}\left[\frac{c(\xi)}{|\xi|^4}\right](x_0)|_{|\xi'|=1}&=&\pi^+_{\xi_n}\left[\frac{c(\xi')+\xi_nc(dx_n)}{(1+\xi_n^2)^2}\right]\nonumber\\
&=&\frac{1}{2\pi i}{\rm lim}_{u\rightarrow
0^-}\int_{\Gamma^+}\frac{\frac{c(\xi')+\eta_nc(dx_n)}{(\eta_n+i)^2(\xi_n+iu-\eta_n)}}
{(\eta_n-i)^2}d\eta_n\nonumber\\
%&=&\left[\frac{c(\xi')+\eta_nc(dx_n)}{(\eta_n+i)^2(\xi_n-\eta_n)}\right]^{(1)}|_{\eta_n=i}\nonumber\\
&=&-\frac{(i\xi_n+2)c(\xi')+ic(dx_n)}{4(\xi_n-i)^2}.
\end{eqnarray}
Similarly we have,
\begin{eqnarray}
\pi^+_{\xi_n}\left[\frac{i\partial_{x_n}c(\xi')}{|\xi|^2}\right](x_0)|_{|\xi'|=1}=\frac{\partial_{x_n}[c(\xi')](x_0)}{2(\xi_n-i)}.
\end{eqnarray}
By (3.25), then\\
\begin{eqnarray}\pi^+_{\xi_n}\partial_{x_n}\sigma_{-1}(\widehat{D}^{-1})|_{|\xi'|=1}
=\frac{\partial_{x_n}[c(\xi')](x_0)}{2(\xi_n-i)}+ih'(0)
\left[\frac{(i\xi_n+2)c(\xi')+ic(dx_n)}{4(\xi_n-i)^2}\right].
\end{eqnarray}
\noindent By the relation of the Clifford action and ${\rm tr}{AB}={\rm tr }{BA}$, we have the equalities:\\
\begin{eqnarray*}{\rm tr}[c(\xi')c(dx_n)]=0;~~{\rm tr}[c(dx_n)^2]=-16;~~{\rm tr}[c(\xi')^2](x_0)|_{|\xi'|=1}=-16;
\end{eqnarray*}
\begin{eqnarray}
{\rm tr}[\partial_{x_n}c(\xi')c(dx_n)]=0;~~{\rm tr}[\partial_{x_n}c(\xi')c(\xi')](x_0)|_{|\xi'|=1}=-8h'(0);~~{\rm tr}[\bar{c}(\widetilde{e_i})\bar{c}(\widetilde{e_j})c(\widetilde{e_k})c(\widetilde{e_l})]=0(i\neq j).
\end{eqnarray}
By (3.27) and a direct computation, we have
\begin{eqnarray}
&&h'(0){\rm tr}\bigg[\frac{(i\xi_n+2)c(\xi')+ic(dx_n)}{4(\xi_n-i)^2}\times
\bigg(\frac{6\xi_nc(dx_n)+2c(\xi')}{(1+\xi_n^2)^2}
-\frac{8\xi_n^2[c(\xi')+\xi_nc(dx_n)]}{(1+\xi_n^2)^3}\bigg)
\bigg](x_0)|_{|\xi'|=1}\nonumber\\
&=&-16h'(0)\frac{-2i\xi_n^2-\xi_n+i}{(\xi_n-i)^4(\xi_n+i)^3}.
\end{eqnarray}
Similarly, we have
\begin{eqnarray}
&&-i{\rm
tr}\bigg[\bigg(\frac{\partial_{x_n}[c(\xi')](x_0)}{2(\xi_n-i)}\bigg)
\times\bigg(\frac{6\xi_nc(dx_n)+2c(\xi')}{(1+\xi_n^2)^2}-\frac{8\xi_n^2[c(\xi')+\xi_nc(dx_n)]}
{(1+\xi_n^2)^3}\bigg)\bigg](x_0)|_{|\xi'|=1}\nonumber\\
&=&-8ih'(0)\frac{3\xi_n^2-1}{(\xi_n-i)^4(\xi_n+i)^3}.
\end{eqnarray}
Then\\
\begin{eqnarray*}
{\rm case~a)~II)}&=&-\int_{|\xi'|=1}\int^{+\infty}_{-\infty}\frac{4ih'(0)(\xi_n-i)^2}
{(\xi_n-i)^4(\xi_n+i)^3}d\xi_n\sigma(\xi')dx'\\
&=&-4ih'(0)\Omega_3\int_{\Gamma^+}\frac{1}{(\xi_n-i)^2(\xi_n+i)^3}d\xi_ndx'\\
&=&-4ih'(0)\Omega_32\pi i[\frac{1}{(\xi_n+i)^3}]^{(1)}|_{\xi_n=i}dx'\\
&=&-\frac{3}{2}\pi h'(0)\Omega_3dx'.
\end{eqnarray*}
where ${\rm \Omega_{3}}$ is the canonical volume of $S^{3}.$\\

\noindent  {\bf case a)~III)}~$r=-1,~l=-1,~j=|\alpha|=0,~k=1$\\

\noindent By (3.14), we get
\begin{equation}
{\rm case~a)~III)}=-\frac{1}{2}\int_{|\xi'|=1}\int^{+\infty}_{-\infty}
{\rm trace} [\partial_{\xi_n}\pi^+_{\xi_n}\sigma_{-1}(\widehat{D}^{-1})\times
\partial_{\xi_n}\partial_{x_n}\sigma_{-1}((\widehat{D}^*)^{-1})](x_0)d\xi_n\sigma(\xi')dx'.
\end{equation}
\noindent By Lemma 3.7, we have\\
\begin{eqnarray}
\partial_{\xi_n}\partial_{x_n}\sigma_{-1}((\widehat{D}^*)^{-1})(x_0)|_{|\xi'|=1}
&=&-ih'(0)\left[\frac{c(dx_n)}{|\xi|^4}-4\xi_n\frac{c(\xi')
+\xi_nc(dx_n)}{|\xi|^6}\right]-\frac{2\xi_ni\partial_{x_n}c(\xi')(x_0)}{|\xi|^4};
\end{eqnarray}
\begin{eqnarray}
\partial_{\xi_n}\pi^+_{\xi_n}\sigma_{-1}(\widehat{D}^{-1})(x_0)|_{|\xi'|=1}&=&-\frac{c(\xi')+ic(dx_n)}{2(\xi_n-i)^2}.
\end{eqnarray}
Similar to {\rm case~a)~II)}, we have\\
\begin{eqnarray}
&&{\rm tr}\left\{\frac{c(\xi')+ic(dx_n)}{2(\xi_n-i)^2}\times
ih'(0)\left[\frac{c(dx_n)}{|\xi|^4}-4\xi_n\frac{c(\xi')+\xi_nc(dx_n)}{|\xi|^6}\right]\right\}
=8h'(0)\frac{i-3\xi_n}{(\xi_n-i)^4(\xi_n+i)^3}
\end{eqnarray}
and
\begin{eqnarray}
{\rm tr}\left[\frac{c(\xi')+ic(dx_n)}{2(\xi_n-i)^2}\times
\frac{2\xi_ni\partial_{x_n}c(\xi')(x_0)}{|\xi|^4}\right]
=\frac{-8ih'(0)\xi_n}{(\xi_n-i)^4(\xi_n+i)^2}.
\end{eqnarray}
So we have
\begin{eqnarray}
{\rm case~a)~III)}&=&-\int_{|\xi'|=1}\int^{+\infty}_{-\infty}\frac{h'(0)4(i-3\xi_n)}
{(\xi_n-i)^4(\xi_n+i)^3}d\xi_n\sigma(\xi')dx'
-\int_{|\xi'|=1}\int^{+\infty}_{-\infty}\frac{h'(0)4i\xi_n}
{(\xi_n-i)^4(\xi_n+i)^2}d\xi_n\sigma(\xi')dx'\nonumber\\
&=&-h'(0)\Omega_3\frac{2\pi i}{3!}[\frac{4(i-3\xi_n)}{(\xi_n+i)^3}]^{(3)}|_{\xi_n=i}dx'+h'(0)\Omega_3\frac{2\pi i}{3!}[\frac{4i\xi_n}{(\xi_n+i)^2}]^{(3)}|_{\xi_n=i}dx'\nonumber\\
&=&\frac{3}{2}\pi h'(0)\Omega_3dx'.
\end{eqnarray}

\noindent  {\bf case b)}~$r=-2,~l=-1,~k=j=|\alpha|=0$\\

\noindent By (3.14), we get
\begin{eqnarray}
{\rm case~b)}&=&-i\int_{|\xi'|=1}\int^{+\infty}_{-\infty}{\rm trace} [\pi^+_{\xi_n}\sigma_{-2}(\widehat{D}^{-1})\times
\partial_{\xi_n}\sigma_{-1}((\widehat{D}^*)^{-1})](x_0)d\xi_n\sigma(\xi')dx'.
\end{eqnarray}
 By Lemma 3.7 we have\\
\begin{eqnarray}
\sigma_{-2}(\widehat{D}^{-1})(x_0)=\frac{c(\xi)\sigma_{0}(\widehat{D})(x_0)c(\xi)}{|\xi|^4}+\frac{c(\xi)}{|\xi|^6}c(dx_n)
[\partial_{x_n}[c(\xi')](x_0)|\xi|^2-c(\xi)h'(0)|\xi|^2_{\partial
M}].
\end{eqnarray}
where
\begin{eqnarray}
\sigma_{0}(\widehat{D})(x_0)&=&\frac{1}{4}\sum_{s,t,i}\omega_{s,t}(\widetilde{e_i})
(x_{0})c(\widetilde{e_i})\bar{c}(\widetilde{e_s})\bar{c}(\widetilde{e_t})
-\frac{1}{4}\sum_{s,t,i}\omega_{s,t}(\widetilde{e_i})(x_{0})
c(\widetilde{e_i})c(\widetilde{e_s})c(\widetilde{e_t}))\nonumber\\
&&+\bar{c}(\theta)+c(\theta').
\end{eqnarray}
We denote
\begin{eqnarray}
b_0^{1}(x_0)&=&\frac{1}{4}\sum_{s,t,i}\omega_{s,t}(\widetilde{e_i})
(x_{0})c(\widetilde{e_i})\bar{c}(\widetilde{e_s})\bar{c}(\widetilde{e_t})\nonumber\\
b_0^{2}(x_0)&=&-\frac{1}{4}\sum_{s,t,i}\omega_{s,t}(\widetilde{e_i})(x_{0})
c(\widetilde{e_i})c(\widetilde{e_s})c(\widetilde{e_t})).
\end{eqnarray}
Then
\begin{eqnarray}
&&\pi^+_{\xi_n}\sigma_{-2}(\widehat{D}^{-1}(x_0))|_{|\xi'|=1}\nonumber\\
&=&\pi^+_{\xi_n}\Big[\frac{c(\xi)b_0^{1}(x_0)c(\xi)}{(1+\xi_n^2)^2}\Big]+\pi^+_{\xi_n}
\Big[\frac{c(\xi)(\bar{c}(\theta)+c(\theta'))(x_0)c(\xi)}{(1+\xi_n^2)^2}\Big]
\nonumber\\
&&+\pi^+_{\xi_n}\Big[\frac{c(\xi)b_0^{2}(x_0)c(\xi)+c(\xi)c(dx_n)\partial_{x_n}[c(\xi')](x_0)}{(1+\xi_n^2)^2}-h'(0)\frac{c(\xi)c(dx_n)c(\xi)}{(1+\xi_n^{2})^3}\Big].
\end{eqnarray}
By direct calculation we have

\begin{eqnarray}
\pi^+_{\xi_n}\Big[\frac{c(\xi)b_0^{1}(x_0)c(\xi)}{(1+\xi_n^2)^2}\Big]&=&\pi^+_{\xi_n}\Big[\frac{c(\xi')p_0^{1}(x_0)c(\xi')}{(1+\xi_n^2)^2}\Big]
+\pi^+_{\xi_n}\Big[ \frac{\xi_nc(\xi')b_0^{1}(x_0)c(dx_{n})}{(1+\xi_n^2)^2}\Big]\nonumber\\
&&+\pi^+_{\xi_n}\Big[\frac{\xi_nc(dx_{n})b_0^{1}(x_0)c(\xi')}{(1+\xi_n^2)^2}\Big]
+\pi^+_{\xi_n}\Big[\frac{\xi_n^{2}c(dx_{n})b_0^{1}(x_0)c(dx_{n})}{(1+\xi_n^2)^2}\Big]\nonumber\\
&=&-\frac{c(\xi')b_0^{1}(x_0)c(\xi')(2+i\xi_{n})}{4(\xi_{n}-i)^{2}}
+\frac{ic(\xi')b_0^{1}(x_0)c(dx_{n})}{4(\xi_{n}-i)^{2}}\nonumber\\
&&+\frac{ic(dx_{n})b_0^{1}(x_0)c(\xi')}{4(\xi_{n}-i)^{2}}
+\frac{-i\xi_{n}c(dx_{n})b_0^{1}(x_0)c(dx_{n})}{4(\xi_{n}-i)^{2}}.
\end{eqnarray}
Since
\begin{eqnarray}
c(dx_n)b_0^{1}(x_0)
&=&-\frac{1}{4}h'(0)\sum^{n-1}_{i=1}c(\widetilde{e_i})
\bar{c}(\widetilde{e_i})c(\widetilde{e_n})\bar{c}(\widetilde{e_n}),
\end{eqnarray}
then by the relation of the Clifford action and ${\rm tr}{AB}={\rm tr }{BA}$,  we have the equalities:\\
\begin{eqnarray}
{\rm tr}[c(\widetilde{e_i})
\bar{c}(\widetilde{e_i})c(\widetilde{e_n})
\bar{c}(\widetilde{e_n})]&=&0~~(i<n);~~
{\rm tr}[b_0^{1}c(dx_n)]=0;~~{\rm tr }[\bar{c}(\theta)c(dx_{n})]=0;\nonumber\\
~~{\rm tr }[c(\theta')c(dx_{n})]&=&-16g(\theta',dx_n);
~~{\rm tr}[\bar {c}(\xi')\bar {c}(dx_n)]=0.
\end{eqnarray}
Since
\begin{eqnarray}
\partial_{\xi_n}\sigma_{-1}((\widehat{D}^*)^{-1})=\partial_{\xi_n}q_{-1}(x_0)|_{|\xi'|=1}=i\left[\frac{c(dx_n)}{1+\xi_n^2}-\frac{2\xi_nc(\xi')+2\xi_n^2c(dx_n)}{(1+\xi_n^2)^2}\right],
\end{eqnarray}
By (3.43) and (3.46), we have
\begin{eqnarray}
&&{\rm tr }[\pi^+_{\xi_n}\Big[\frac{c(\xi)b_0^{1}(x_0)c(\xi)}{(1+\xi_n^2)^2}\Big]
\times\partial_{\xi_n}\sigma_{-1}((\widehat{D}^*)^{-1})(x_0)]|_{|\xi'|=1}\nonumber\\
&=&\frac{1}{2(1+\xi_n^2)^2}{\rm tr }[c(\xi')b_0^{1}(x_0)]
+\frac{i}{2(1+\xi_n^2)^2}{\rm tr }[c(dx_n)b_0^{1}(x_0)]\nonumber\\
&=&\frac{1}{2(1+\xi_n^2)^2}{\rm tr }[c(\xi')b_0^{1}(x_0)].
\end{eqnarray}
We note that $i<n,~\int_{|\xi'|=1}\{\xi_{i_{1}}\xi_{i_{2}}\cdots\xi_{i_{2d+1}}\}\sigma(\xi')=0$,
so ${\rm tr }[c(\xi')b_0^{1}(x_0)]$ has no contribution for computing {\rm case~b)}.

By direct calculation we have
\begin{eqnarray}
\pi^+_{\xi_n}\Big[\frac{c(\xi)b_0^{2}(x_0)c(\xi)+c(\xi)c(dx_n)\partial_{x_n}[c(\xi')](x_0)}{(1+\xi_n^2)^2}\Big]-h'(0)\pi^+_{\xi_n}\Big[\frac{c(\xi)c(dx_n)c(\xi)}{(1+\xi_n)^3}\Big]:= B_1-B_2,
\end{eqnarray}
where
\begin{eqnarray}
B_1&=&\frac{-1}{4(\xi_n-i)^2}[(2+i\xi_n)c(\xi')b_0^{2}(x_0)c(\xi')+i\xi_nc(dx_n)b_0^{2}(x_0)c(dx_n)\nonumber\\
&&+(2+i\xi_n)c(\xi')c(dx_n)\partial_{x_n}c(\xi')+ic(dx_n)b_0^{2}(x_0)c(\xi')
+ic(\xi')b_0^{2}(x_0)c(dx_n)-i\partial_{x_n}c(\xi')]
\end{eqnarray}
and
\begin{eqnarray}
B_2&=&\frac{h'(0)}{2}\left[\frac{c(dx_n)}{4i(\xi_n-i)}+\frac{c(dx_n)-ic(\xi')}{8(\xi_n-i)^2}
+\frac{3\xi_n-7i}{8(\xi_n-i)^3}[ic(\xi')-c(dx_n)]\right].
\end{eqnarray}
By (3.46) and (3.50), we have\\
\begin{eqnarray}{\rm tr }[B_2\times\partial_{\xi_n}\sigma_{-1}((\widehat{D}^*)^{-1})]|_{|\xi'|=1}
&=&\frac{i}{2}h'(0)\frac{-i\xi_n^2-\xi_n+4i}{4(\xi_n-i)^3(\xi_n+i)^2}{\rm tr}[ \texttt{id}]\nonumber\\
&=&8ih'(0)\frac{-i\xi_n^2-\xi_n+4i}{4(\xi_n-i)^3(\xi_n+i)^2}.
\end{eqnarray}
By (3.46) and (3.49), we have
\begin{eqnarray}{\rm tr }[B_1\times\partial_{\xi_n}\sigma_{-1}((\widehat{D}^*)^{-1})]|_{|\xi'|=1}=
\frac{-8ic_0}{(1+\xi_n^2)^2}+2h'(0)\frac{\xi_n^2-i\xi_n-2}{(\xi_n-i)(1+\xi_n^2)^2},
\end{eqnarray}
where $b_0^{2}=c_0c(dx_n)$ and $c_0=-\frac{3}{4}h'(0)$.\\

By (3.52) and (3.51), we have
\begin{eqnarray}
&&-i\int_{|\xi'|=1}\int^{+\infty}_{-\infty}{\rm trace} [(B_1-B_2)\times
\partial_{\xi_n}\sigma_{-1}((\widehat{D}^*)^{-1})](x_0)d\xi_n\sigma(\xi')dx'\nonumber\\
&=&-\Omega_3\int_{\Gamma^+}\frac{8c_0(\xi_n-i)+ih'(0)}{(\xi_n-i)^3(\xi_n+i)^2}d\xi_ndx'\nonumber\\
&=&\frac{9}{2}\pi h'(0)\Omega_3dx'.
\end{eqnarray}
Similar to (3.47), we have
\begin{eqnarray}
&&{\rm tr }[\pi^+_{\xi_n}\Big[\frac{c(\xi)\bar{c}(\theta)(x_0)c(\xi)}{(1+\xi_n^2)^2}\Big]
\times\partial_{\xi_n}\sigma_{-1}((\widehat{D}^*)^{-1})(x_0)]|_{|\xi'|=1}\nonumber\\
&=&\frac{1}{2(1+\xi_n^2)^2}{\rm tr }[c(\xi')\bar{c}(\theta)(x_0)]
+\frac{i}{2(1+\xi_n^2)^2}{\rm tr }[c(dx_n)\bar{c}(\theta)(x_0)]\nonumber\\
&=&\frac{1}{2(1+\xi_n^2)^2}{\rm tr }[c(\xi')\bar{c}(\theta)(x_0)].
\end{eqnarray}
Similar to (3.49), we have
\begin{eqnarray}
{\rm tr }[\pi^+_{\xi_n}\Big[\frac{c(\xi)c(\theta')(x_0)c(\xi)}{(1+\xi_n^2)^2}\Big]
\times\partial_{\xi_n}\sigma_{-1}((\widehat{D}^*)^{-1})(x_0)]|_{|\xi'|=1}
&=&\frac{i}{2(1+\xi_n^2)^2}{\rm tr }[c(dx_n)c(\theta')(x_0)].
\end{eqnarray}
By (3.54) and (3.55), we have
\begin{eqnarray}
&&-i\int_{|\xi'|=1}\int^{+\infty}_{-\infty}{\rm trace} [\pi^+_{\xi_n}
\Big[\frac{c(\xi)(\bar{c}(\theta)+c(\theta'))c(\xi)}{(1+\xi_n^2)^2}\Big]\times
\partial_{\xi_n}\sigma_{-1}((\widehat{D}^*)^{-1})](x_0)d\xi_n\sigma(\xi')dx'\nonumber\\
&=&\frac{\pi}{4}{\rm tr}[c(dx_n)c(\theta')]\Omega_3dx'\nonumber\\
&=&-4\pi g(\theta',dx_n)\Omega_3dx'.
\end{eqnarray}
By (3.53) and (3.56), we have\\
\begin{eqnarray}
{\rm case~b)}=\frac{9}{2}\pi h'(0)\Omega_3dx'-4\pi g(\theta',dx_n)\Omega_3dx'.
\end{eqnarray}

\noindent {\bf  case c)}~$r=-1,~l=-2,~k=j=|\alpha|=0$\\
By (3.14), we get
\begin{equation}
{\rm case~ c)}=-i\int_{|\xi'|=1}\int^{+\infty}_{-\infty}{\rm trace} [\pi^+_{\xi_n}\sigma_{-1}(\widehat{D}^{-1})\times
\partial_{\xi_n}\sigma_{-2}((\widehat{D}^*)^{-1})](x_0)d\xi_n\sigma(\xi')dx'.
\end{equation}
By (3.2) and (3.3), Lemma 3.7, we have
\begin{equation}
\pi^+_{\xi_n}\sigma_{-1}(\widehat{D}^{-1})|_{|\xi'|=1}=\frac{c(\xi')+ic(dx_n)}{2(\xi_n-i)}.
\end{equation}

Since
\begin{equation}
\sigma_{-2}((\widehat{D}^*)^{-1})(x_0)=\frac{c(\xi)\sigma_{0}(\widehat{D}^*)(x_0)c(\xi)}{|\xi|^4}+\frac{c(\xi)}{|\xi|^6}c(dx_n)
\bigg[\partial_{x_n}[c(\xi')](x_0)|\xi|^2-c(\xi)h'(0)|\xi|^2_{\partial
M}\bigg],
\end{equation}
where
\begin{eqnarray}
\sigma_{0}(\widehat{D}^*)(x_0)&=&\frac{1}{4}\sum_{s,t,i}\omega_{s,t}(\widetilde{e_i})(x_{0})c(\widetilde{e_i})\bar{c}(\widetilde{e_s})\bar{c}(\widetilde{e_t})
-\frac{1}{4}\sum_{s,t,i}\omega_{s,t}(\widetilde{e_i})(x_{0})c(\widetilde{e_i})
c(\widetilde{e_s})c(\widetilde{e_t}))+(\bar{c}(\theta)-c(\theta'))(x_{0})\nonumber\\
&=&b_0^{1}(x_0)+b_0^{2}(x_0)+(\bar{c}(\theta)-c(\theta'))(x_{0}),
\end{eqnarray}
then
\begin{eqnarray}
\partial_{\xi_n}\sigma_{-2}((\widehat{D}^*)^{-1})(x_0)|_{|\xi'|=1}&=&
\partial_{\xi_n}\bigg\{\frac{c(\xi)[b_0^{1}(x_0)+b_0^{2}(x_0)
+(\bar{c}(\theta)-c(\theta'))(x_{0})]c(\xi)}{|\xi|^4}\nonumber\\
&&+\frac{c(\xi)}{|\xi|^6}c(dx_n)[\partial_{x_n}[c(\xi')](x_0)|\xi|^2-c(\xi)h'(0)]\bigg\}\nonumber\\
&=&\partial_{\xi_n}\bigg\{\frac{c(\xi)b_0^{1}(x_0)]c(\xi)}{|\xi|^4}+\frac{c(\xi)}{|\xi|^6}c(dx_n)[\partial_{x_n}[c(\xi')](x_0)|\xi|^2-c(\xi)h'(0)]\bigg\}\nonumber\\
&&+\partial_{\xi_n}\frac{c(\xi)b_0^{2}(x_0)c(\xi)}{|\xi|^4}
+\partial_{\xi_n}\frac{c(\xi)(\bar{c}(\theta)-c(\theta'))(x_{0})c(\xi)}{|\xi|^4}.
\end{eqnarray}
By direct calculation we have
\begin{eqnarray}
\partial_{\xi_n}\frac{c(\xi)b_0^{1}(x_0)c(\xi)}{|\xi|^4}=\frac{c(dx_n)b_0^{1}(x_0)c(\xi)}{|\xi|^4}
+\frac{c(\xi)b_0^{1}(x_0)c(dx_n)}{|\xi|^4}
-\frac{4\xi_n c(\xi)b_0^{1}(x_0)c(\xi)}{|\xi|^6};
\end{eqnarray}
\begin{eqnarray}
\partial_{\xi_n}\frac{c(\xi)(\bar{c}(\theta)-c(\theta'))(x_{0})c(\xi)}{|\xi|^4}
&=&\frac{c(dx_n)(\bar{c}(\theta)-c(\theta'))(x_{0})c(\xi)}{|\xi|^4}
+\frac{c(\xi)(\bar{c}(\theta)-c(\theta'))(x_{0})c(dx_n)}{|\xi|^4}\nonumber\\
&-&\frac{4\xi_n c(\xi)(\bar{c}(\theta)-c(\theta'))(x_{0})c(\xi)}{|\xi|^4}.
\end{eqnarray}
We denote $$q_{-2}^{1}=\frac{c(\xi)b_0^{2}(x_0)c(\xi)}{|\xi|^4}+\frac{c(\xi)}{|\xi|^6}c(dx_n)[\partial_{x_n}[c(\xi')](x_0)|\xi|^2-c(\xi)h'(0)],$$ then
\begin{eqnarray}
\partial_{\xi_n}(q_{-2}^{1})&=&\frac{1}{(1+\xi_n^2)^3}\bigg[(2\xi_n-2\xi_n^3)c(dx_n)b_0^{2}c(dx_n)
+(1-3\xi_n^2)c(dx_n)b_0^{2}c(\xi')\nonumber\\
&&+(1-3\xi_n^2)c(\xi')b_0^{2}c(dx_n)
-4\xi_nc(\xi')b_0^{2}c(\xi')
+(3\xi_n^2-1)\partial_{x_n}c(\xi')\nonumber\\
&&-4\xi_nc(\xi')c(dx_n)\partial_{x_n}c(\xi')
+2h'(0)c(\xi')+2h'(0)\xi_nc(dx_n)\bigg]\nonumber\\
&&+6\xi_nh'(0)\frac{c(\xi)c(dx_n)c(\xi)}{(1+\xi^2_n)^4};
\end{eqnarray}
By (3.59) and (3.63), we have
\begin{eqnarray}
&&{\rm tr}[\pi^+_{\xi_n}\sigma_{-1}(\widehat{D}^{-1})\times
\partial_{\xi_n}\frac{c(\xi)b_0^{1}c(\xi)}
{|\xi|^4}](x_0)|_{|\xi'|=1}\nonumber\\
&=&\frac{-1}{(\xi-i)(\xi+i)^3}{\rm tr}[c(\xi')b_0^{1}(x_0)]
+\frac{i}{(\xi-i)(\xi+i)^3}{\rm tr}[c(dx_n)b_0^{1}(x_0)].
\end{eqnarray}
By (3.45), we have
\begin{eqnarray}
{\rm tr}[\pi^+_{\xi_n}\sigma_{-1}(\widehat{D}^{-1})\times
\partial_{\xi_n}\frac{c(\xi)b_0^{1}c(\xi)}
{|\xi|^4}](x_0)|_{|\xi'|=1}
=\frac{-1}{(\xi-i)(\xi+i)^3}{\rm tr}[c(\xi')b_0^{1}(x_0)].
\end{eqnarray}
We note that $i<n,~\int_{|\xi'|=1}\{\xi_{i_{1}}\xi_{i_{2}}\cdots\xi_{i_{2d+1}}\}\sigma(\xi')=0$,
so ${\rm tr }[c(\xi')b_0^{1}(x_0)]$ has no contribution for computing {\rm case~c)}.

By (3.59) and (3.65), we have
\begin{eqnarray}
{\rm tr}[\pi^+_{\xi_n}\sigma_{-1}(\widehat{D}^{-1})\times
\partial_{\xi_n}(q^1_{-2})](x_0)|_{|\xi'|=1}
=\frac{12h'(0)(i\xi^2_n+\xi_n-2i)}{(\xi-i)^3(\xi+i)^3}
+\frac{48h'(0)i\xi_n}{(\xi-i)^3(\xi+i)^4},
\end{eqnarray}
then
\begin{eqnarray}
-i\Omega_3\int_{\Gamma_+}[\frac{12h'(0)(i\xi_n^2+\xi_n-2i)}
{(\xi_n-i)^3(\xi_n+i)^3}+\frac{48h'(0)i\xi_n}{(\xi_n-i)^3(\xi_n+i)^4}]d\xi_ndx'=
-\frac{9}{2}\pi h'(0)\Omega_3dx'.
\end{eqnarray}
By (3.59) and (3.64), we have
\begin{eqnarray}
&&{\rm tr}[\pi^+_{\xi_n}\sigma_{-1}(\widehat{D}^{-1})\times
\partial_{\xi_n}\frac{c(\xi)(\bar{c}(\theta)-c(\theta'))c(\xi)}
{|\xi|^4}](x_0)|_{|\xi'|=1}\nonumber\\
&=&\frac{-1}{(\xi-i)(\xi+i)^3}{\rm tr}[c(\xi')(\bar{c}(\theta)-c(\theta'))(x_0)]
+\frac{i}{(\xi-i)(\xi+i)^3}{\rm tr}[c(dx_n)(\bar{c}(\theta)-c(\theta'))(x_0)].
\end{eqnarray}

By $\int_{|\xi'|=1}\{\xi_{1}\cdot\cdot\cdot\xi_{2d+1}\}\sigma(\xi')=0$ and (3.47), we have
\begin{eqnarray}
&&-i\int_{|\xi'|=1}\int^{+\infty}_{-\infty}{\rm tr}[\pi^+_{\xi_n}\sigma_{-1}(\widehat{D}^{-1})\times
\partial_{\xi_n}\frac{c(\xi)(\bar{c}(\theta)-c(\theta'))c(\xi)}
{|\xi|^4}](x_0)d\xi_n\sigma(\xi')dx'\nonumber\\
&=&-i\int_{|\xi'|=1}\int^{+\infty}_{-\infty}\frac{i}{(\xi-i)(\xi+i)^3}{\rm tr}[c(dx_n)(\bar{c}(\theta)-c(\theta'))](x_0)d\xi_n\sigma(\xi')dx'\nonumber\\
&=&-\frac{\pi}{4}{\rm tr}[c(dx_n)(\bar{c}(\theta)-c(\theta'))]\Omega_3dx'\nonumber\\
&=&-4\pi g(\theta',dx_n)\Omega_3dx'.
\end{eqnarray}
So we have
\begin{eqnarray}
{\rm case~ c)}=-\frac{9}{2}\pi h'(0)\Omega_3dx'-4\pi g(\theta',dx_n)\Omega_3dx'.
\end{eqnarray}
Since $\Phi$ is the sum of the cases a), b) and c), so $\Phi=-8\pi g(\theta',dx_n)\Omega_3dx'$.

\begin{thm}
Let $M$ be a $4$-dimensional oriented
compact manifolds with the boundary $\partial M$ and the metric
$g^M$ as above, $\widehat{D}$ and $\widehat{D}^*$ be modified Novikov operators on $\widehat{M}$, then
\begin{eqnarray}
&&\widetilde{{\rm Wres}}[\pi^+\widehat{D}^{-1}\circ\pi^+(\widehat{D}^*)^{-1}]\nonumber\\
&=&32\pi^2\int_{M}
\bigg[16g(\widetilde{e_{j}},\nabla^{TM}_{\widetilde{e_{j}}}\theta')-\frac{4}{3}s
-16|\theta|^2+32|\theta'|^2\bigg]d{\rm Vol_{M}}-8\pi\int_{\partial M}g(\theta',dx_n)\Omega_3dx'.
\end{eqnarray}
where $s$ is the scalar curvature.
\end{thm}
On the other hand, we also prove the Kastler-Kalau-Walze type theorem for $4$-dimensional manifolds with boundary associated to $\widehat{D}^2$.
By (3.6) and (3.7), we will compute
\begin{equation}
\widetilde{{\rm Wres}}[\pi^+\widehat{D}^{-1}\circ\pi^+\widehat{D}^{-1}]=\int_M\int_{|\xi|=1}{\rm
trace}_{\wedge^*T^*M}[\sigma_{-4}(\widehat{D}^{-2})]\sigma(\xi)dx+\int_{\partial M}\widehat{\Phi},
\end{equation}
where
\begin{eqnarray}
\widehat{\Phi} &=&\int_{|\xi'|=1}\int^{+\infty}_{-\infty}\sum^{\infty}_{j, k=0}\sum\frac{(-i)^{|\alpha|+j+k+1}}{\alpha!(j+k+1)!}
\times {\rm trace}_{\wedge^*T^*M}[\partial^j_{x_n}\partial^\alpha_{\xi'}\partial^k_{\xi_n}\sigma^+_{r}(\widehat{D}^{-1})(x',0,\xi',\xi_n)
\nonumber\\
&&\times\partial^\alpha_{x'}\partial^{j+1}_{\xi_n}\partial^k_{x_n}\sigma_{l}(\widehat{D}^{-1})(x',0,\xi',\xi_n)]d\xi_n\sigma(\xi')dx',
\end{eqnarray}
and the sum is taken over $r+l-k-j-|\alpha|=-3,~~r\leq -1,l\leq-1$.\\

Locally we can use Theorem 2.2 (2.29) to compute the interior of $\widetilde{{\rm Wres}}[\pi^+\widehat{D}^{-1}\circ\pi^+\widehat{D}^{-1}]$, we have
\begin{eqnarray}
&&\int_M\int_{|\xi|=1}{\rm
trace}_{\wedge^*T^*M}[\sigma_{-4}(\widehat{D}^{-2})]\sigma(\xi)dx=32\pi^2\int_{M}
\bigg[-\frac{4}{3}s
-16|\theta|^2\bigg]d{\rm Vol_{M}}.
\end{eqnarray}

So we only need to compute $\int_{\partial M} \widehat{\Phi}$. From the remark above, now we can compute $\widehat{\Phi}$ (see formula (3.75) for the definition of $\widehat{\Phi}$). We use ${\rm tr}$ as shorthand of ${\rm trace}$. Since $n=4$, then ${\rm tr}_{\wedge^*T^*M}[{\rm \texttt{id}}]={\rm dim}(\wedge^*(4))=16$, since the sum is taken over $
r+l-k-j-|\alpha|=-3,~~r\leq -1,l\leq-1,$ then we have the following five cases:

~\\
\noindent  {\bf case a)~I)}~$r=-1,~l=-1,~k=j=0,~|\alpha|=1$\\

\noindent By (3.75), we get
\begin{equation}
{\rm case~a)~I)}=-\int_{|\xi'|=1}\int^{+\infty}_{-\infty}\sum_{|\alpha|=1}
 {\rm tr}[\partial^\alpha_{\xi'}\pi^+_{\xi_n}\sigma_{-1}(\widehat{D}^{-1})\times
 \partial^\alpha_{x'}\partial_{\xi_n}
 \sigma_{-1}(\widehat{D}^{-1})](x_0)d\xi_n\sigma(\xi')dx'.
\end{equation}
\noindent  {\bf case a)~II)}~$r=-1,~l=-1,~k=|\alpha|=0,~j=1$\\

\noindent By (3.75), we get
\begin{equation}
{\rm case~ a)~II)}=-\frac{1}{2}\int_{|\xi'|=1}\int^{+\infty}_{-\infty} {\rm
trace} [\partial_{x_n}\pi^+_{\xi_n}\sigma_{-1}(\widehat{D}^{-1})\times
\partial_{\xi_n}^2\sigma_{-1}(\widehat{D}^{-1})](x_0)d\xi_n\sigma(\xi')dx'.
\end{equation}

\noindent  {\bf case a)~III)}~$r=-1,~l=-1,~j=|\alpha|=0,~k=1$\\

\noindent By (3.75), we get
\begin{equation}
{\rm case~ a)~III)}=-\frac{1}{2}\int_{|\xi'|=1}\int^{+\infty}_{-\infty}
{\rm trace} [\partial_{\xi_n}\pi^+_{\xi_n}\sigma_{-1}(\widehat{D}^{-1})\times
\partial_{\xi_n}\partial_{x_n}\sigma_{-1}(\widehat{D}^{-1})](x_0)d\xi_n\sigma(\xi')dx'.
\end{equation}
By Lemma 3.7, we have $\sigma_{-1}(\widehat{D}^{-1})=\sigma_{-1}((\widehat{D}^*)^{-1})$.
By (3.21)-(3.37), so {\bf case a)} vanishes.\\

\noindent  {\bf case b)}~$r=-2,~l=-1,~k=j=|\alpha|=0$\\

\noindent By (3.75), we get
\begin{eqnarray}
{\rm case~ b)}&=&-i\int_{|\xi'|=1}\int^{+\infty}_{-\infty}{\rm trace} [\pi^+_{\xi_n}\sigma_{-2}(\widehat{D}^{-1})\times
\partial_{\xi_n}\sigma_{-1}(\widehat{D}^{-1})](x_0)d\xi_n\sigma(\xi')dx'.
\end{eqnarray}
By Lemma 3.7, we have $\sigma_{-1}(\widehat{D}^{-1})=\sigma_{-1}((\widehat{D}^*)^{-1})$.
By (3.38)-(3.57), we have\\
\begin{eqnarray}
{\rm case~ b)}=\frac{9}{2}\pi h'(0)\Omega_3dx'-4\pi g(\theta',dx_n)\Omega_3dx',
\end{eqnarray}
where ${\rm \Omega_{4}}$ is the canonical volume of $S^{4}.$\\
\noindent {\bf  case c)}~$r=-1,~l=-2,~k=j=|\alpha|=0$\\
By (3.75), we get
\begin{equation}
{\rm case~ c)}=-i\int_{|\xi'|=1}\int^{+\infty}_{-\infty}{\rm trace} [\pi^+_{\xi_n}\sigma_{-1}(\widehat{D}^{-1})\times
\partial_{\xi_n}\sigma_{-2}(\widehat{D}^{-1})](x_0)d\xi_n\sigma(\xi')dx'.
\end{equation}
By (3.2) and (3.3), Lemma 3.7, we have
\begin{equation}
\pi^+_{\xi_n}\sigma_{-1}(\widehat{D}^{-1})|_{|\xi'|=1}=\frac{c(\xi')+ic(dx_n)}{2(\xi_n-i)}.
\end{equation}

Since
\begin{equation}
\sigma_{-2}(\widehat{D}^{-1})(x_0)=\frac{c(\xi)\sigma_{0}(\widehat{D})(x_0)c(\xi)}{|\xi|^4}+\frac{c(\xi)}{|\xi|^6}c(dx_n)
[\partial_{x_n}[c(\xi')](x_0)|\xi|^2-c(\xi)h'(0)|\xi|^2_{\partial
M}],
\end{equation}
where
\begin{eqnarray}
\sigma_{0}(\widehat{D})(x_0)&=&\frac{1}{4}\sum_{s,t,i}\omega_{s,t}(\widetilde{e_i})(x_{0})c(\widetilde{e_i})\bar{c}(\widetilde{e_s})\bar{c}(\widetilde{e_t})
-\frac{1}{4}\sum_{s,t,i}\omega_{s,t}(\widetilde{e_i})(x_{0})c(\widetilde{e_i})
c(\widetilde{e_s})c(\widetilde{e_t}))+(\bar{c}(\theta)+c(\theta'))(x_{0})\nonumber\\
&=&b_0^{1}(x_0)+b_0^{2}(x_0)+(\bar{c}(\theta)+c(\theta'))(x_{0}),
\end{eqnarray}
then
\begin{eqnarray}
\partial_{\xi_n}\sigma_{-2}(\widehat{D}^{-1})(x_0)|_{|\xi'|=1}&=&
\partial_{\xi_n}\bigg\{\frac{c(\xi)[b_0^{1}(x_0)+b_0^{2}(x_0)
+(\bar{c}(\theta)+c(\theta'))(x_{0})]c(\xi)}{|\xi|^4}\nonumber\\
&&+\frac{c(\xi)}{|\xi|^6}c(dx_n)[\partial_{x_n}[c(\xi')](x_0)|\xi|^2-c(\xi)h'(0)]\bigg\}\nonumber\\
&=&\partial_{\xi_n}\bigg\{\frac{c(\xi)b_0^{1}(x_0)]c(\xi)}{|\xi|^4}+\frac{c(\xi)}{|\xi|^6}c(dx_n)[\partial_{x_n}[c(\xi')](x_0)|\xi|^2-c(\xi)h'(0)]\bigg\}\nonumber\\
&&+\partial_{\xi_n}\frac{c(\xi)b_0^{2}(x_0)c(\xi)}{|\xi|^4}
+\partial_{\xi_n}\frac{c(\xi)(\bar{c}(\theta)+c(\theta'))(x_{0})c(\xi)}{|\xi|^4}.
\end{eqnarray}
By direct calculation we have
\begin{eqnarray}
\partial_{\xi_n}\frac{c(\xi)b_0^{1}(x_0)c(\xi)}{|\xi|^4}=\frac{c(dx_n)b_0^{1}(x_0)c(\xi)}{|\xi|^4}
+\frac{c(\xi)b_0^{1}(x_0)c(dx_n)}{|\xi|^4}
-\frac{4\xi_n c(\xi)b_0^{1}(x_0)c(\xi)}{|\xi|^6};
\end{eqnarray}
\begin{eqnarray}
\partial_{\xi_n}\frac{c(\xi)(\bar{c}(\theta)+c(\theta'))(x_{0})c(\xi)}{|\xi|^4}
&&=\frac{c(dx_n)(\bar{c}(\theta)+c(\theta'))(x_{0})c(\xi)}{|\xi|^4}
+\frac{c(\xi)(\bar{c}(\theta)+c(\theta'))(x_{0})c(dx_n)}{|\xi|^4}\nonumber\\
&&-\frac{4\xi_n c(\xi)(\bar{c}(\theta)+c(\theta'))(x_{0})c(\xi)}{|\xi|^4}.
\end{eqnarray}
We denote $$q_{-2}^{1}=\frac{c(\xi)b_0^{2}(x_0)c(\xi)}{|\xi|^4}+\frac{c(\xi)}{|\xi|^6}c(dx_n)[\partial_{x_n}[c(\xi')](x_0)|\xi|^2-c(\xi)h'(0)],$$ then
\begin{eqnarray}
\partial_{\xi_n}(q_{-2}^{1})&=&\frac{1}{(1+\xi_n^2)^3}\bigg[(2\xi_n-2\xi_n^3)c(dx_n)b_0^{2}c(dx_n)
+(1-3\xi_n^2)c(dx_n)b_0^{2}c(\xi')\nonumber\\
&&+ (1-3\xi_n^2)c(\xi')b_0^{2}c(dx_n)
-4\xi_nc(\xi')b_0^{2}c(\xi')
+(3\xi_n^2-1)\partial_{x_n}c(\xi')\nonumber\\
&&-4\xi_nc(\xi')c(dx_n)\partial_{x_n}c(\xi')
+2h'(0)c(\xi')+2h'(0)\xi_nc(dx_n)\bigg]\nonumber\\
&&+6\xi_nh'(0)\frac{c(\xi)c(dx_n)c(\xi)}{(1+\xi^2_n)^4}.
\end{eqnarray}
By (3.83) and (3.87), we have
\begin{eqnarray}
&&{\rm tr}[\pi^+_{\xi_n}\sigma_{-1}(\widehat{D}^{-1})\times
\partial_{\xi_n}\frac{c(\xi)b_0^{1}c(\xi)}
{|\xi|^4}](x_0)|_{|\xi'|=1}\nonumber\\
&=&\frac{-1}{(\xi-i)(\xi+i)^3}{\rm tr}[c(\xi')b_0^{1}(x_0)]
+\frac{i}{(\xi-i)(\xi+i)^3}{\rm tr}[c(dx_n)b_0^{1}(x_0)].
\end{eqnarray}
By (3.45), we have
\begin{eqnarray}
{\rm tr}[\pi^+_{\xi_n}\sigma_{-1}(\widehat{D}^{-1})\times
\partial_{\xi_n}\frac{c(\xi)b_0^{1}c(\xi)}
{|\xi|^4}](x_0)|_{|\xi'|=1}
=\frac{-1}{(\xi-i)(\xi+i)^3}{\rm tr}[c(\xi')b_0^{1}(x_0)].
\end{eqnarray}
We note that $i<n,~\int_{|\xi'|=1}\{\xi_{i_{1}}\xi_{i_{2}}\cdots\xi_{i_{2d+1}}\}\sigma(\xi')=0$,
so ${\rm tr }[c(\xi')b_0^{1}(x_0)]$ has no contribution for computing case c).

By (3.83) and (3.89), we have
\begin{eqnarray}
{\rm tr}[\pi^+_{\xi_n}\sigma_{-1}(\widehat{D}^{-1})\times
\partial_{\xi_n}(q^1_{-2})](x_0)|_{|\xi'|=1}
=\frac{12h'(0)(i\xi^2_n+\xi_n-2i)}{(\xi-i)^3(\xi+i)^3}
+\frac{48h'(0)i\xi_n}{(\xi-i)^3(\xi+i)^4},
\end{eqnarray}
then
\begin{eqnarray}
-i\Omega_3\int_{\Gamma_+}[\frac{12h'(0)(i\xi_n^2+\xi_n-2i)}
{(\xi_n-i)^3(\xi_n+i)^3}+\frac{48h'(0)i\xi_n}{(\xi_n-i)^3(\xi_n+i)^4}]d\xi_ndx'=
-\frac{9}{2}\pi h'(0)\Omega_3dx'.
\end{eqnarray}
By (3.83) and (3.88), we have
\begin{eqnarray}
&&{\rm tr}[\pi^+_{\xi_n}\sigma_{-1}(\widehat{D}^{-1})\times
\partial_{\xi_n}\frac{c(\xi)(\bar{c}(\theta)+c(\theta'))c(\xi)}
{|\xi|^4}](x_0)|_{|\xi'|=1}\nonumber\\
&=&\frac{-1}{(\xi-i)(\xi+i)^3}{\rm tr}[c(\xi')(\bar{c}(\theta)+c(\theta'))(x_0)]
+\frac{i}{(\xi-i)(\xi+i)^3}{\rm tr}[c(dx_n)(\bar{c}(\theta)+c(\theta'))(x_0)].
\end{eqnarray}

By $\int_{|\xi'|=1}\{\xi_{1}\cdot\cdot\cdot\xi_{2d+1}\}\sigma(\xi')=0$ and (3.45), we have
\begin{eqnarray}
&&-i\int_{|\xi'|=1}\int^{+\infty}_{-\infty}{\rm tr}[\pi^+_{\xi_n}\sigma_{-1}(\widehat{D}^{-1})\times
\partial_{\xi_n}\frac{c(\xi)(\bar{c}(\theta)+c(\theta'))c(\xi)}
{|\xi|^4}](x_0)d\xi_n\sigma(\xi')dx'\nonumber\\
&=&-i\int_{|\xi'|=1}\int^{+\infty}_{-\infty}\frac{i}{(\xi-i)(\xi+i)^3}{\rm tr}[c(dx_n)(\bar{c}(\theta)+c(\theta'))](x_0)d\xi_n\sigma(\xi')dx'\nonumber\\
&=&-\frac{\pi}{4}{\rm tr}[c(dx_n)(\bar{c}(\theta)+c(\theta'))]\Omega_3dx'\nonumber\\
&=&4\pi g(\theta',dx_n)\Omega_3dx'.
\end{eqnarray}
So we have
\begin{eqnarray}
{\rm case~ c)}=-\frac{9}{2}\pi h'(0)\Omega_3dx'+4\pi g(\theta',dx_n)\Omega_3dx'.
\end{eqnarray}
Since $\widehat{\Phi}$ is the sum of the cases a), b) and c), so $\widehat{\Phi}=0$.

\begin{thm}
Let $M$ be a $4$-dimensional oriented
compact manifold with the boundary $\partial M$ and the metric
$g^M$ as above, $\widehat{D}$ be modified Novikov operator on $\widehat{M}$, then
\begin{eqnarray}
&&\widetilde{{\rm Wres}}[\pi^+\widehat{D}^{-1}\circ\pi^+\widehat{D}^{-1}]=32\pi^2\int_{M}\bigg(-\frac{4}{3}s
-16|\theta|^2\bigg)d{\rm Vol_{M}}.
\end{eqnarray}
where $s$ is the scalar curvature.
\end{thm}

%%%%%%%%%第4部分%%%%%%%%
\section{A Kastler-Kalau-Walze type theorem for $6$-dimensional manifolds with boundary }
In this section, we prove the Kastler-Kalau-Walze type theorems for $6$-dimensional manifolds with boundary. An application of (2.1.4) in \cite{Wa5} shows that

\begin{equation}
\widetilde{{\rm Wres}}[\pi^+\widehat{D}^{-1}\circ\pi^+(\widehat{D}^{*}\widehat{D}
      \widehat{D}^{*})^{-1}]=\int_M\int_{|\xi|=1}{\rm
trace}_{\wedge ^*T^*M}[\sigma_{-4}((\widehat{D}^*\widehat{D})^{-2})]\sigma(\xi)dx+\int_{\partial M}\Psi,
\end{equation}
where
\begin{eqnarray}
\Psi &=&\int_{|\xi'|=1}\int^{+\infty}_{-\infty}\sum^{\infty}_{j, k=0}\sum\frac{(-i)^{|\alpha|+j+k+1}}{\alpha!(j+k+1)!}
\times {\rm trace}_{\wedge ^*T^*M}[\partial^j_{x_n}\partial^\alpha_{\xi'}\partial^k_{\xi_n}\sigma^+_{r}(\widehat{D}^{-1})(x',0,\xi',\xi_n)
\nonumber\\
&&\times\partial^\alpha_{x'}\partial^{j+1}_{\xi_n}\partial^k_{x_n}\sigma_{l}
((\widehat{D}^{*}\widehat{D}
      \widehat{D}^{*})^{-1})(x',0,\xi',\xi_n)]d\xi_n\sigma(\xi')dx',
\end{eqnarray}
and the sum is taken over $r+\ell-k-j-|\alpha|-1=-6, \ r\leq-1, \ell\leq -3$.\\
Locally we can use Theorem 2.2 (2.28) to compute the interior term of (4.1), we have
\begin{eqnarray}
&&\int_M\int_{|\xi|=1}{\rm
trace}_{\wedge^*T^*M}[\sigma_{-4}((\widehat{D}^*\widehat{D})^{-2})]\sigma(\xi)dx\nonumber\\
&=&128\pi^3\int_{M}\bigg[64g(\widetilde{e_{j}},\nabla^{TM}
_{\widetilde{e_{j}}}\theta')-\frac{16}{3}s
-64|\theta|^2+256|\theta'|^2\bigg]d{\rm Vol_{M}}.
\end{eqnarray}

So we only need to compute $\int_{\partial M} \Psi$. Let us now turn to compute the specification of
$\widehat{D}^*\widehat{D}\widehat{D}^*$.
\begin{eqnarray}
\widehat{D}^*\widehat{D}\widehat{D}^*
&=&\sum^{n}_{i=1}c(\widetilde{e_{i}})\langle \widetilde{e_{i}},dx_{l}\rangle(-g^{ij}\partial_{l}\partial_{i}\partial_{j})
+\sum^{n}_{i=1}c(\widetilde{e_{i}})\langle \widetilde{e_{i}},dx_{l}\rangle \bigg\{-(\partial_{l}g^{ij})\partial_{i}\partial_{j}-g^{ij}\bigg(4(\sigma_{i}
+a_{i})\partial_{j}-2\Gamma^{k}_{ij}\partial_{k}\bigg)\partial_{l}\bigg\} \nonumber\\
&&+\sum^{n}_{i=1}c(\widetilde{e_{i}})\langle \widetilde{e_{i}},dx_{l}\rangle \bigg\{-2(\partial_{l}g^{ij})(\sigma_{i}+a_i)\partial_{j}+g^{ij}
(\partial_{l}\Gamma^{k}_{ij})\partial_{k}-2g^{ij}[(\partial_{l}\sigma_{i})
+(\partial_{l}a_i)]\partial_{j}
         +(\partial_{l}g^{ij})\Gamma^{k}_{ij}\partial_{k}\nonumber\\
         &&+\sum_{j,k}\Big[\partial_{l}\Big(c(\theta')c(\widetilde{e_{j}})
         -c(\widetilde{e_{j}}) c(\theta')\Big)\Big]\langle \widetilde{e_{j}},dx^{k}\rangle\partial_{k}
         +\sum_{j,k}\Big(c(\theta')c(\widetilde{e_{j}})-c(\widetilde{e_{j}})c(\theta')\Big)\Big[\partial_{l}\langle \widetilde{e_{j}},dx^{k}\rangle\Big]\partial_{k} \bigg\}\nonumber\\
         &&+\sum^{n}_{i=1}c(\widetilde{e_{i}})\langle \widetilde{e_{i}},dx_{l}\rangle\partial_{l}\bigg\{-g^{ij}\Big[
         (\partial_{i}\sigma_{j})+(\partial_{i}a_{j})+\sigma_{i}\sigma_{j}+\sigma_{i}a_{j}
         +a_{i}\sigma_{j}+a_{i}a_{j}-\Gamma_{i,j}^{k}\sigma_{k}-\Gamma_{i,j}^{k}a_{k}\nonumber\\
         &&+\sum_{i,j}g^{i,j}\Big[c(\theta')c(\partial_{i})\sigma_{i}
         +c(\theta')c(\partial_{i})a_{i}-c(\partial_{i})\partial_{i}(c(\theta'))
         -c(\partial_{i})\sigma_{i}c(\theta')-c(\partial_{i})a_{i}c(\theta')\Big]+\frac{1}{4}s\nonumber\\
         &&-\frac{1}{8}\sum_{ijkl}R_{ijkl}\bar{c}(\widetilde{e_i})\bar{c}(\widetilde{e_j})
         c(\widetilde{e_k})c(\widetilde{e_l})
         +\sum_{i}c(\widetilde{e_i})\bar{c}(\nabla_{\widetilde{e_i}}^{TM}\theta)
         +|\theta|^{2}+|\theta'|^2-\bar{c}(\theta)c(\theta')+c(\theta')\bar{c}(\theta)\bigg\}\nonumber\\
         &&+\Big[(\sigma_{i}+a_{i})+(\bar{c}(\theta)-c(\theta'))\Big](-g^{ij}\partial_{i}\partial_{j})
         +\sum^{n}_{i=1}c(\widetilde{e_{i}})\langle \widetilde{e_{i}},dx_{l}\rangle \bigg\{2\sum_{j,k}\Big[c(\theta')c(\widetilde{e_{j}})-c(\widetilde{e_{j}})c(\theta')\Big]\nonumber\\
         &&\times\langle \widetilde{e_{i}},dx_{k}\rangle\bigg\}\partial_{l}\partial_{k}
         +\Big[(\sigma_{i}+a_{i})+(\bar{c}(\theta)-c(\theta'))\Big]
         \bigg\{-\sum_{i,j}g^{i,j}\Big[2\sigma_{i}\partial_{j}+2a_{i}\partial_{j}
         -\Gamma_{i,j}^{k}\partial_{k}+(\partial_{i}\sigma_{j})\nonumber\\
         &&+(\partial_{i}a_{j})+\sigma_{i}\sigma_{j}+\sigma_{i}a_{j}+a_{i}\sigma_{j}+a_{i}a_{j} -\Gamma_{i,j}^{k}\sigma_{k}-\Gamma_{i,j}^{k}a_{k}\Big]-\sum_{i,j}g^{i,j}\Big[c(\partial_{i})c(\theta')
         -c(\theta')c(\partial_{i})\Big]\partial_{j}\nonumber\\
         &&+\sum_{i,j}g^{i,j}\Big[c(\theta')c(\partial_{i})\sigma_{i}+c(\theta')c(\partial_{i})a_{i}
         -c(\partial_{i})\partial_{i}(c(\theta'))-c(\partial_{i})\sigma_{i}c(\theta')
         -c(\partial_{i})a_{i}c(\theta')\Big]+\frac{1}{4}s+|\theta'|^2\nonumber\\
         &&-\frac{1}{8}\sum_{ijkl}R_{ijkl}\bar{c}(\widetilde{e_i})\bar{c}(\widetilde{e_j})
         c(\widetilde{e_k})c(\widetilde{e_l})
         +\sum_{i}c(\widetilde{e_i})\bar{c}(\nabla_{\widetilde{e_i}}^{TM}\theta)
         +|\theta|^{2}-\bar{c}(\theta)c(\theta')+c(\theta')\bar{c}(\theta)\bigg\}.
\end{eqnarray}
Then, we obtain
\begin{lem} The following identities hold:
\begin{eqnarray}
\sigma_2(\widehat{D}^*\widehat{D}\widehat{D}^*)&=&\sum_{i,j,l}c(dx_{l})\partial_{l}(g^{i,j})\xi_{i}\xi_{j} +c(\xi)(4\sigma^k+4a^k-2\Gamma^k)\xi_{k}-2[c(\xi)c(\theta')c(\xi)+|\xi|^2c(\theta')]\nonumber\\
&&+\frac{1}{4}|\xi|^2\sum_{s,t,l}\omega_{s,t}
(\widetilde{e_l})[c(\widetilde{e_l})\bar{c}(\widetilde{e_s})\bar{c}(\widetilde{e_t})
-c(\widetilde{e_l})c(\widetilde{e_s})c(\widetilde{e_t})]
+|\xi|^2(\bar{c}(\theta)-c(\theta'));\nonumber\\
\sigma_{3}(\widehat{D}^*\widehat{D}\widehat{D}^*)
&=&ic(\xi)|\xi|^{2}.
\end{eqnarray}
\end{lem}

Write
\begin{eqnarray}
\sigma(\widehat{D}^*\widehat{D}\widehat{D}^*)&=&p_3+p_2+p_1+p_0;
~\sigma((\widehat{D}^*\widehat{D}\widehat{D}^*)^{-1})=\sum^{\infty}_{j=3}q_{-j}.
\end{eqnarray}

By the composition formula of pseudodifferential operators, we have

\begin{eqnarray}
1=\sigma((\widehat{D}^*\widehat{D}\widehat{D}^*)\circ (\widehat{D}^*\widehat{D}\widehat{D}^*)^{-1})&=&
\sum_{\alpha}\frac{1}{\alpha!}\partial^{\alpha}_{\xi}
[\sigma(\widehat{D}^*\widehat{D}\widehat{D}^*)]D^{\alpha}_{x}
[(\widehat{D}^*\widehat{D}\widehat{D}^*)^{-1}] \nonumber\\
&=&(p_3+p_2+p_1+p_0)(q_{-3}+q_{-4}+q_{-5}+\cdots) \nonumber\\
&+&\sum_j(\partial_{\xi_j}p_3+\partial_{\xi_j}p_2+\partial_{\xi_j}p_1+\partial_{\xi_j}p_0)
(D_{x_j}q_{-3}+D_{x_j}q_{-4}+D_{x_j}q_{-5}+\cdots) \nonumber\\
&=&p_3q_{-3}+(p_3q_{-4}+p_2q_{-3}+\sum_j\partial_{\xi_j}p_3D_{x_j}q_{-3})+\cdots,
\end{eqnarray}
by (4.7), we have

\begin{equation}
q_{-3}=p_3^{-1};~q_{-4}=-p_3^{-1}[p_2p_3^{-1}+\sum_j\partial_{\xi_j}p_3D_{x_j}(p_3^{-1})].
\end{equation}
By Lemma 4.1, we have some symbols of operators.
\begin{lem} The following identities hold:
\begin{eqnarray}
\sigma_{-3}((\widehat{D}^*\widehat{D}\widehat{D}^*)^{-1})&=&\frac{ic(\xi)}{|\xi|^{4}};\nonumber\\
\sigma_{-4}((\widehat{D}^*\widehat{D}\widehat{D}^*)^{-1})&=&
\frac{c(\xi)\sigma_{2}(\widehat{D}^*\widehat{D}\widehat{D}^*)c(\xi)}{|\xi|^8}
+\frac{ic(\xi)}{|\xi|^8}\Big(|\xi|^4c(dx_n)\partial_{x_n}c(\xi')
-2h'(0)c(dx_n)c(\xi)\nonumber\\
&&+2\xi_{n}c(\xi)\partial_{x_n}c(\xi')+4\xi_{n}h'(0)\Big).
\end{eqnarray}
\end{lem}

From the remark above, now we can compute $\Psi$ (see formula (4.2) for the definition of $\Psi$). We use ${\rm tr}$ as shorthand of ${\rm trace}$. Since $n=6$, then ${\rm tr}_{\wedge ^*T^*M}[\texttt{id}]=64$.
Since the sum is taken over $r+\ell-k-j-|\alpha|-1=-6, \ r\leq-1, \ell\leq -3$, then we have the $\int_{\partial_{M}}\Psi$
is the sum of the following five cases:

~\\
\noindent  {\bf case (a)~(I)}~$r=-1, l=-3, j=k=0, |\alpha|=1$.\\
By (4.2), we get
 \begin{equation}
{\rm case~(a)~(I)}=-\int_{|\xi'|=1}\int^{+\infty}_{-\infty}\sum_{|\alpha|=1}{\rm trace}
\Big[\partial^{\alpha}_{\xi'}\pi^{+}_{\xi_{n}}\sigma_{-1}(\widehat{D}^{-1})
      \times\partial^{\alpha}_{x'}\partial_{\xi_{n}}\sigma_{-3}((\widehat{D}^{*}\widehat{D}
      \widehat{D}^{*})^{-1})\Big](x_0)d\xi_n\sigma(\xi')dx'.
\end{equation}
By Lemma 4.2, for $i<n$, we have
 \begin{equation}
 \partial_{x_{i}}\sigma_{-3}((\widehat{D}^{*}\widehat{D}\widehat{D}^{*})^{-1})(x_0)=
      \partial_{x_{i}}\Big[\frac{ic(\xi)}{|\xi|^{4}}\Big](x_{0})
      =i\partial_{x_{i}}[c(\xi)]|\xi|^{-4}(x_{0})-2ic(\xi)\partial_{x_{i}}[|\xi|^{2}]|\xi|^{-6}(x_{0})=0.
\end{equation}
 so {\rm case~(a)~(I)} vanishes.
~\\

\noindent  {\bf case (a)~(II)}~$r=-1, l=-3, |\alpha|=k=0, j=1$.\\
By (4.2), we have
  \begin{equation}
{\rm case~(a)~(II)}=-\frac{1}{2}\int_{|\xi'|=1}\int^{+\infty}_{-\infty} {\rm
trace} \Big[\partial_{x_{n}}\pi^{+}_{\xi_{n}}\sigma_{-1}(\widehat{D}^{-1})
      \times\partial^{2}_{\xi_{n}}\sigma_{-3}((\widehat{D}^{*}\widehat{D}\widehat{D}^{*})^{-1})\Big](x_0)d\xi_n\sigma(\xi')dx'.
\end{equation}
By Lemma 4.2 and direct calculations, we have\\
\begin{equation}
\partial^{2}_{\xi_{n}}\sigma_{-3}((\widehat{D}^{*}\widehat{D}\widehat{D}^{*})^{-1})=i\bigg[\frac{(20\xi^{2}_{n}-4)c(\xi')+
12(\xi^{3}_{n}-\xi_{n})c(dx_{n})}{(1+\xi_{n}^{2})^{4}}\bigg].
\end{equation}
Since $n=6$, ${\rm tr}[-\texttt{id}]=-64$. By the relation of the Clifford action and ${\rm tr}AB={\rm tr}BA$,  then
\begin{eqnarray}
&&{\rm tr}[c(\xi')c(dx_{n})]=0; \ {\rm tr}[c(dx_{n})^{2}]=-64;\
{\rm tr}[c(\xi')^{2}](x_{0})|_{|\xi'|=1}=-64;\nonumber\\
&&{\rm tr}[\partial_{x_{n}}[c(\xi')]c(\texttt{d}x_{n})]=0; \
{\rm tr}[\partial_{x_{n}}c(\xi')c(\xi')](x_{0})|_{|\xi'|=1}=-32h'(0).
\end{eqnarray}
By (3.28), (4.13) and (4.14), we get
\begin{equation}
{\rm
trace} \Big[\partial_{x_{n}}\pi^{+}_{\xi_{n}}\sigma_{-1}(\widehat{D}^{-1})
      \times\partial^{2}_{\xi_{n}}\sigma_{-3}((\widehat{D}^{*}\widehat{D}\widehat{D}^{*})^{-1})\Big](x_0)
=64 h'(0)\frac{-1-3\xi_{n}i+5\xi^{2}_{n}+3i\xi^{3}_{n}}{(\xi_{n}-i)^{6}(\xi_{n}+i)^{4}}.
\end{equation}
Then we obtain

\begin{eqnarray}
{\rm case~(a)~(II)}&=&-\frac{1}{2}\int_{|\xi'|=1}\int^{+\infty}_{-\infty} h'(0)dimF\frac{-8-24\xi_{n}i+40\xi^{2}_{n}+24i\xi^{3}_{n}}{(\xi_{n}-i)^{6}(\xi_{n}+i)^{4}}d\xi_n\sigma(\xi')dx'\nonumber\\
     &=&8h'(0)\Omega_{4}\int_{\Gamma^{+}}\frac{4+12\xi_{n}i-20\xi^{2}_{n}-122i\xi^{3}_{n}}{(\xi_{n}-i)^{6}(\xi_{n}+i)^{4}}d\xi_{n}dx'\nonumber\\
     &=&h'(0)\Omega_{4}\frac{\pi i}{5!}\Big[\frac{8+24\xi_{n}i-40\xi^{2}_{n}-24i\xi^{3}_{n}}{(\xi_{n}+i)^{4}}\Big]^{(5)}|_{\xi_{n}=i}dx'\nonumber\\
     &=&-\frac{15}{2}\pi h'(0)\Omega_{4}dx',
\end{eqnarray}
where ${\rm \Omega_{4}}$ is the canonical volume of $S^{4}.$\\

\noindent  {\bf case (a)~(III)}~$r=-1,l=-3,|\alpha|=j=0,k=1$.\\
By (4.2), we have
 \begin{equation}
{\rm case~ (a)~(III)}=-\frac{1}{2}\int_{|\xi'|=1}\int^{+\infty}_{-\infty}{\rm trace} \Big[\partial_{\xi_{n}}\pi^{+}_{\xi_{n}}\sigma_{-1}(\widehat{D}^{-1})
      \times\partial_{\xi_{n}}\partial_{x_{n}}\sigma_{-3}((\widehat{D}^{*}\widehat{D}\widehat{D}^{*})^{-1})\Big](x_0)d\xi_n\sigma(\xi')dx'.
\end{equation}
By Lemma 4.2 and direct calculations, we have\\
\begin{equation}
\partial_{\xi_{n}}\partial_{x_{n}}\sigma_{-3}((\widehat{D}^{*}\widehat{D}\widehat{D}^{*})^{-1})=-\frac{4 i\xi_{n}\partial_{x_{n}}c(\xi')(x_{0})}{(1+\xi_{n}^{2})^{3}}
      +i\frac{12h'(0)\xi_{n}c(\xi')}{(1+\xi_{n}^{2})^{4}}
      -i\frac{(2-10\xi^{2}_{n})h'(0)c(dx_{n})}{(1+\xi_{n}^{2})^{4}}.
\end{equation}
Combining (3.34) and (4.18), we have
\begin{equation}
{\rm trace} \Big[\partial_{\xi_{n}}\pi^{+}_{\xi_{n}}\sigma_{-1}(\widehat{D}^{-1})
      \times\partial_{\xi_{n}}\partial_{x_{n}}\sigma_{-3}((\widehat{D}^{*}\widehat{D}\widehat{D}^{*})^{-1})\Big](x_{0})|_{|\xi'|=1}
=8h'(0)\frac{8i-32\xi_{n}-8i\xi^{2}_{n}}{(\xi_{n}-i)^{5}(\xi+i)^{4}}.
\end{equation}
Then
\begin{eqnarray}
{\rm case~(a)~III)}&=&-\frac{1}{2}\int_{|\xi'|=1}\int^{+\infty}_{-\infty} 8h'(0)\frac{8i-32\xi_{n}-8i\xi^{2}_{n}}{(\xi_{n}-i)^{5}(\xi+i)^{4}}d\xi_n\sigma(\xi')dx'\nonumber\\
     &=&-\frac{1}{2}h'(0)8\Omega_{4}\int_{\Gamma^{+}}\frac{8i-32\xi_{n}-8i\xi^{2}_{n}}{(\xi_{n}-i)^{5}(\xi+i)^{4}}d\xi_{n}dx'\nonumber\\
     &=&-8h'(0)\Omega_{4}\frac{\pi i}{4!}\Big[\frac{8i-32\xi_{n}-8i\xi^{2}_{n}}{(\xi+i)^{4}}\Big]^{(4)}|_{\xi_{n}=i}dx'\nonumber\\
     &=&\frac{25}{2}\pi h'(0)\Omega_{4}dx'.
\end{eqnarray}

\noindent  {\bf case (b)}~$r=-1,l=-4,|\alpha|=j=k=0$.\\
By (4.2), we have
 \begin{eqnarray}
{\rm case~ (b)}&=&-i\int_{|\xi'|=1}\int^{+\infty}_{-\infty}{\rm trace} \Big[\pi^{+}_{\xi_{n}}\sigma_{-1}(\widehat{D}^{-1})
      \times\partial_{\xi_{n}}\sigma_{-4}((\widehat{D}^{*}\widehat{D}
      \widehat{D}^{*})^{-1})\Big](x_0)d\xi_n\sigma(\xi')dx'\nonumber\\
&=&i\int_{|\xi'|=1}\int^{+\infty}_{-\infty}{\rm trace} [\partial_{\xi_n}\pi^+_{\xi_n}\sigma_{-1}(\widehat{D}^{-1})\times
\sigma_{-4}((\widehat{D}^{*}\widehat{D}
      \widehat{D}^{*})^{-1})](x_0)d\xi_n\sigma(\xi')dx'.
\end{eqnarray}

In the normal coordinate, $g^{ij}(x_{0})=\delta^{j}_{i}$ and $\partial_{x_{j}}(g^{\alpha\beta})(x_{0})=0$, if $j<n$; $\partial_{x_{j}}(g^{\alpha\beta})(x_{0})=h'(0)\delta^{\alpha}_{\beta}$, if $j=n$.
So by Lemma A.2 in \cite{Wa3}, we have $\Gamma^{n}(x_{0})=\frac{5}{2}h'(0)$ and $\Gamma^{k}(x_{0})=0$ for $k<n$. By the definition of $\delta^{k}$ and Lemma 2.3 in \cite{Wa3}, we have $\delta^{n}(x_{0})=0$ and $\delta^{k}=\frac{1}{4}h'(0)c(\widetilde{e_{k}})c(\widetilde{e_{n}})$ for $k<n$. By Lemma 4.2, we obtain

\begin{eqnarray}
\sigma_{-4}((\widehat{D}^{*}\widehat{D}\widehat{D}^{*})^{-1})(x_{0})|_{|\xi'|=1}&=&
\frac{c(\xi)\sigma_{2}((\widehat{D}^{*}\widehat{D}\widehat{D}^{*})^{-1})
(x_{0})|_{|\xi'|=1}c(\xi)}{|\xi|^8}
-\frac{c(\xi)}{|\xi|^4}\sum_j\partial_{\xi_j}\big(c(\xi)|\xi|^2\big)
D_{x_j}\big(\frac{ic(\xi)}{|\xi|^4}\big)\nonumber\\
&=&\frac{1}{|\xi|^8}c(\xi)\Big(\frac{1}{2}h'(0)c(\xi)\sum_{k<n}\xi_k
c(\widetilde{e_k})c(\widetilde{e_n})-\frac{1}{2}h'(0)c(\xi)\sum_{k<n}\xi_k
\bar{c}(\widetilde{e_k})\bar{c}(\widetilde{e_n})\nonumber\\
&&-\frac{5}{2}h'(0)\xi_nc(\xi)-\frac{1}{4}h'(0)|\xi|^2c(dx_n)
-2[c(\xi)c(\theta')c(\xi)+|\xi|^2c(\theta')]\nonumber\\
&&+|\xi|^2(\widehat{c}(\theta)-c(\theta'))\Big)c(\xi)
+\frac{ic(\xi)}{|\xi|^8}\Big(|\xi|^4c(dx_n)\partial_{x_n}c(\xi')
-2h'(0)c(dx_n)c(\xi)\nonumber\\
&&+2\xi_{n}c(\xi)\partial_{x_n}c(\xi')+4\xi_{n}h'(0)\Big).
\end{eqnarray}
By (3.34) and (4.22), we have
\begin{eqnarray}
&&{\rm tr} [\partial_{\xi_n}\pi^+_{\xi_n}\sigma_{-1}(\widehat{D}^{-1})\times
\sigma_{-4}(\widehat{D}^{*}\widehat{D}\widehat{D}^{*})^{-1}](x_0)|_{|\xi'|=1} \nonumber\\
%&&={\rm tr} \Big[-\frac{c(\xi')+ic(dx_n)}{2(\xi_n-i)^2} \times
 %\frac{1}{(1+\xi_{n}^{2})^{4}}\Big\{\Big( (3+4i)\xi_{n}+3\xi_{n}^{3}\Big)h'(0)c(\xi') \nonumber\\
 %&& +\big(\frac{3}{4}-2i+(6i+2)\xi_{n}^{2}+\frac{9}{4}\xi_{n}^{4}\big)h'(0)c(dx_n) \nonumber\\
%&& +\big(-\frac{1}{2}\xi_{n}-\frac{1}{2}\xi_{n}^{3}\big)h'(0)\bar{c}(\xi')\bar{c}(dx_n)c(dx_n) \nonumber\\
%&& +\big(-\frac{1}{2}-\frac{1}{2}\xi_{n}^{2}\big)h'(0)\bar{c}(\xi')\bar{c}(dx_n)c(\xi') \nonumber\\
%&& +\big(1+2\xi_{n}^{2}+\xi_{n}^{4}\big)T\bar{c}(v)+\big(-3i\xi_{n}-4i\xi_{n}^{3}-i\xi_{n}^{5}\big)\partial_{x_n}c(\xi')\big\}\nonumber\\
&&=\frac{1}{2(\xi_{n}-i)^{2}(1+\xi_{n}^{2})^{4}}\big(\frac{3}{4}i+2+(3+4i)\xi_{n}+(-6+2i)\xi_{n}^{2}+3\xi_{n}^{3}+\frac{9i}{4}\xi_{n}^{4}\big)h'(0){\rm tr}
[id]\nonumber\\
&&+\frac{1}{2(\xi_{n}-i)^{2}(1+\xi_{n}^{2})^{4}}\big(-1-3i\xi_{n}-2\xi_{n}^{2}-4i\xi_{n}^{3}-\xi_{n}^{4}-i\xi_{n}^{5}\big){\rm tr[c(\xi')\partial_{x_n}c(\xi')]}\nonumber\\
&&-\frac{1}{2(\xi_{n}-i)^{2}(1+\xi_{n}^{2})^{4}}\big(\frac{1}{2}i+\frac{1}{2}\xi_{n}+\frac{1}{2}\xi_{n}^{2}+\frac{1}{2}\xi_{n}^{3}\big){\rm tr}
[c(\xi')\bar{c}(\xi')c(dx_n)\bar{c}(dx_n)]\nonumber\\
&&+\frac{-\xi_ni+3}{2(\xi_{n}-i)^{4}(i+\xi_{n})^{3}}{\rm tr}\big[c(\theta')c(dx_n)\big]-\frac{3\xi_n+i}{2(\xi_{n}-i)^{4}(i+\xi_{n})^{3}}{\rm tr}\big[c(\theta')c(\xi')\big].
\end{eqnarray}
By direct calculation and the relation of the Clifford action and ${\rm tr}{AB}={\rm tr }{BA}$, then we have equalities:
\begin{eqnarray}
&&{\rm tr }[c(\theta')(x_0)c(dx_n)]=-64g(\theta',dx_n);~~{\rm tr }[c(\theta')(x_0)c(\xi')]=-64g(\theta',\xi');\nonumber\\
&&{\rm tr}[c(\widetilde{e_i})
\bar{c}(\widetilde{e_i})c(\widetilde{e_n})
\bar{c}(\widetilde{e_n})]=0~~(i<n).
\end{eqnarray}
Then
\begin{eqnarray}
{\rm tr}
[c(\xi')\bar{c}(\xi')c(dx_n)\bar{c}(dx_n)]&=&
\sum_{i<n,j<n}{\rm tr}[\xi_{i}\xi_{j}c(\widetilde{e_i})\bar{c}
(\widetilde{e_j})c(dx_n)\bar{c}(dx_n)]=0.
\end{eqnarray}
So, we have
\begin{eqnarray}
{\rm case~ (b)}&=&
 ih'(0)\int_{|\xi'|=1}\int^{+\infty}_{-\infty}64\times\frac{\frac{3}{4}i+2+(3+4i)\xi_{n}+(-6+2i)\xi_{n}^{2}+3\xi_{n}^{3}+\frac{9i}{4}\xi_{n}^{4}}{2(\xi_n-i)^5(\xi_n+i)^4}d\xi_n\sigma(\xi')dx'\nonumber\\ &+&ih'(0)\int_{|\xi'|=1}\int^{+\infty}_{-\infty}32\times\frac{1+3i\xi_{n}
 +2\xi_{n}^{2}+4i\xi_{n}^{3}+\xi_{n}^{4}+i\xi_{n}^{5}}{2(\xi_{n}-i)^{2}
 (1+\xi_{n}^{2})^{4}}d\xi_n\sigma(\xi')dx'\nonumber\\
 &+&i\int_{|\xi'|=1}\int^{+\infty}_{-\infty}
 \frac{\xi_n-i-2\xi_n i +1}{2(\xi_{n}-i)^{4}(i+\xi_{n})^{3}}{\rm tr}\big[c(\theta')c(dx_n)\big]d\xi_n\sigma(\xi')dx'\nonumber\\
 &-&i\int_{|\xi'|=1}\int^{+\infty}_{-\infty}
 \frac{3\xi_n+i}{2(\xi_{n}-i)^{4}(i+\xi_{n})^{3}}{\rm tr}\big[c(\theta')c(\xi')\big]d\xi_n\sigma(\xi')dx'\nonumber\\
&=&(-\frac{19}{4}i-15)\pi h'(0)\Omega_4dx'+(-\frac{3}{8}i-\frac{75}{8})\pi h'(0)\Omega_4dx'+120i\pi g(dx_n,\theta')\Omega_4dx'\nonumber\\
&=&(-\frac{41}{8}i-\frac{195}{8})\pi h'(0)\Omega_4dx'+120i\pi g(dx_n,\theta')\Omega_4dx'.
\end{eqnarray}

\noindent {\bf  case (c)}~$r=-2,l=-3,|\alpha|=j=k=0$.\\
By (4.2), we have

\begin{equation}
{\rm case~ (c)}=-i\int_{|\xi'|=1}\int^{+\infty}_{-\infty}{\rm trace} \Big[\pi^{+}_{\xi_{n}}\sigma_{-2}(\widehat{D}^{-1})
      \times\partial_{\xi_{n}}\sigma_{-3}((\widehat{D}^{*}\widehat{D}\widehat{D}^{*})^{-1})\Big](x_0)d\xi_n\sigma(\xi')dx'.
\end{equation}

By Lemma 4.1 and Lemma 4.2, we have
\begin{eqnarray}
\sigma_{-2}(\widehat{D}^{-1})(x_0)&=&\frac{c(\xi)\sigma_{0}(\widehat{D})c(\xi)}{|\xi|^4}(x_0)+\frac{c(\xi)}{|\xi|^6}\sum_jc(dx_j)
\Big[\partial_{x_j}(c(\xi))|\xi|^2-c(\xi)\partial_{x_j}(|\xi|^2)\Big](x_0),
\end{eqnarray}
where
\begin{equation}
\sigma_0(\widehat{D})=\frac{1}{4}\sum_{i,s,t}\omega_{s,t}(\widetilde{e_i})c(\widetilde{e_i})\bar{c}(\widetilde{e_s})\bar{c}(\widetilde{e_t})
-\frac{1}{4}\sum_{i,s,t}\omega_{s,t}(\widetilde{e_i})c(\widetilde{e_i})
c(\widetilde{e_s})c(\widetilde{e_t})+\widehat{c}(\theta)+c(\theta').
\end{equation}

On the other hand,
\begin{equation}
\partial_{\xi_{n}}\sigma_{-3}((\widehat{D}^{*}\widehat{D}\widehat{D}^{*})^{-1})=\frac{-4 i \xi_{n}c(\xi')}{(1+\xi_{n}^{2})^{3}}+\frac{i(1- 3\xi_{n}^{2})c(\texttt{d}x_{n})}
{(1+\xi_{n}^{2})^{3}}.
\end{equation}
By (4.28) and (3.2), (3.3), we have
\begin{eqnarray}
\pi^{+}_{\xi_{n}}\Big(\sigma_{-2}(\widehat{D}^{-1})\Big)(x_{_{0}})|_{|\xi'|=1}
&=&\pi^{+}_{\xi_{n}}\Big[\frac{c(\xi)\sigma_{0}(\widehat{D})(x_{0})c(\xi)
+c(\xi)c(dx_{n})\partial_{x_{n}}[c(\xi')](x_{0})}{(1+\xi^{2}_{n})^{2}}\Big]\nonumber\\
&&-h'(0)\pi^{+}_{\xi_{n}}\Big[\frac{c(\xi)c(dx_{n})c(\xi)}{(1+\xi^{2}_{n})^{3}}\Big].
\end{eqnarray}
We denote
 \begin{eqnarray}
\sigma_{0}(\widehat{D})(x_0)|_{\xi_n=i}=b_0(x_0)=b_0^{1}(x_0)+b_0^{2}(x_0)
+\bar{c}(\theta)+c(\theta').
\end{eqnarray}
Then, we obtain
\begin{eqnarray}
\pi^{+}_{\xi_{n}}\Big(\sigma_{-2}(\widehat{D}^{-1})\Big)(x_{_{0}})|_{|\xi'|=1}
&=&\pi^+_{\xi_n}\Big[\frac{c(\xi)b_0^{2}(x_0)c(\xi)+c(\xi)c(dx_n)
\partial_{x_n}[c(\xi')](x_0)}{(1+\xi_n^2)^2}-h'(0)\frac{c(\xi)c(dx_n)c(\xi)}{(1+\xi_n^{2})^3}\Big]\nonumber\\
&&+\pi^+_{\xi_n}\Big[\frac{c(\xi)[b_0^{1}(x_0)]c(\xi)(x_0)}{(1+\xi_n^2)^2}\Big]
+\pi^+_{\xi_n}\Big[\frac{c(\xi)[\bar{c}(\theta)+c(\theta')]c(\xi)(x_0)}{(1+\xi_n^2)^2}\Big].
\end{eqnarray}

Furthermore,
\begin{eqnarray}
&&\pi^+_{\xi_n}\Big[\frac{c(\xi)[\bar{c}(\theta)+c(\theta')](x_0)c(\xi)}
{(1+\xi_n^2)^2}\Big]\nonumber\\
&=&\pi^+_{\xi_n}\Big[\frac{c(\xi')[\bar{c}(\theta)+c(\theta')](x_0)c(\xi')}
{(1+\xi_n^2)^2}\Big]
+\pi^+_{\xi_n}\Big[ \frac{\xi_nc(\xi')[\bar{c}(\theta)+c(\theta')]
(x_0)c(dx_{n})}{(1+\xi_n^2)^2}\Big]\nonumber\\
&&+\pi^+_{\xi_n}\Big[\frac{\xi_nc(dx_{n})[\bar{c}(\theta)
+c(\theta')](x_0)c(\xi')}{(1+\xi_n^2)^2}\Big]
+\pi^+_{\xi_n}\Big[\frac{\xi_n^{2}c(dx_{n})[\bar{c}(\theta)
+c(\theta')](x_0)c(dx_{n})}{(1+\xi_n^2)^2}\Big]\nonumber\\
&=&-\frac{c(\xi')[\bar{c}(\theta)+c(\theta')](x_0)c(\xi')(2+i\xi_{n})}{4(\xi_{n}-i)^{2}}
+\frac{ic(\xi')[\bar{c}(\theta)+c(\theta')](x_0)c(dx_{n})}{4(\xi_{n}-i)^{2}}\nonumber\\
&&+\frac{ic(dx_{n})[\bar{c}(\theta)+c(\theta')](x_0)c(\xi')}{4(\xi_{n}-i)^{2}}
+\frac{-i\xi_{n}c(dx_{n})[\bar{c}(\theta)+c(\theta')](x_0)c(dx_{n})}{4(\xi_{n}-i)^{2}},
\end{eqnarray}
\begin{eqnarray}
\pi^+_{\xi_n}\Big[\frac{c(\xi)b_0^{1}(x_0)c(\xi)}{(1+\xi_n^2)^2}\Big]&=&\pi^+_{\xi_n}\Big[\frac{c(\xi')p_0^{1}(x_0)c(\xi')}{(1+\xi_n^2)^2}\Big]
+\pi^+_{\xi_n}\Big[ \frac{\xi_nc(\xi')b_0^{1}(x_0)c(dx_{n})}{(1+\xi_n^2)^2}\Big]\nonumber\\
&&+\pi^+_{\xi_n}\Big[\frac{\xi_nc(dx_{n})b_0^{1}(x_0)c(\xi')}{(1+\xi_n^2)^2}\Big]
+\pi^+_{\xi_n}\Big[\frac{\xi_n^{2}c(dx_{n})b_0^{1}(x_0)c(dx_{n})}{(1+\xi_n^2)^2}\Big]\nonumber\\
&=&-\frac{c(\xi')b_0^{1}(x_0)c(\xi')(2+i\xi_{n})}{4(\xi_{n}-i)^{2}}
+\frac{ic(\xi')b_0^{1}(x_0)c(dx_{n})}{4(\xi_{n}-i)^{2}}\nonumber\\
&&+\frac{ic(dx_{n})b_0^{1}(x_0)c(\xi')}{4(\xi_{n}-i)^{2}}
+\frac{-i\xi_{n}c(dx_{n})b_0^{1}(x_0)c(dx_{n})}{4(\xi_{n}-i)^{2}}.
\end{eqnarray}
By the relation of the Clifford action and ${\rm tr}{AB}={\rm tr }{BA}$, then we have equalities:
\begin{eqnarray}
{\rm tr}[b_0^{1}c(dx_n)]=0;~~{\rm tr}[\bar {c}(\xi')\bar {c}(dx_n)]=0;
~~{\rm tr}[\bar{c}(\theta)c(\xi')]=0.
\end{eqnarray}

Then we have
\begin{eqnarray}
{\rm tr }\bigg[\pi^+_{\xi_n}\Big(\frac{c(\xi)b_0^{1}(x_0)c(\xi)}{(1+\xi_n^2)^2}\Big)\times
\partial_{\xi_n}\sigma_{-3}((\widehat{D}^{*}\widehat{D}\widehat{D}^{*})^{-1})(x_0)\bigg]\bigg|_{|\xi'|=1}=\frac{2-8i\xi_n-6\xi_n^2}{4(\xi_n-i)^{2}(1+\xi_n^2)^{3}}{\rm tr }[b_0^{1}(x_0)c(\xi')],
\end{eqnarray}

By direct calculation we have
\begin{eqnarray}
\pi^+_{\xi_n}\Big[\frac{c(\xi)b_0^{2}(x_0)c(\xi)+c(\xi)c(dx_n)\partial_{x_n}(c(\xi'))(x_0)}{(1+\xi_n^2)^2}\Big]-h'(0)\pi^+_{\xi_n}\Big[\frac{c(\xi)c(dx_n)c(\xi)}{(1+\xi_n)^3}\Big]:= B_1-B_2,\nonumber\\
\end{eqnarray}
where
\begin{eqnarray}
B_1&=&\frac{-1}{4(\xi_n-i)^2}\big[(2+i\xi_n)c(\xi')b^2_0c(\xi')+i\xi_nc(dx_n)b^2_0c(dx_n) \nonumber\\
&&+(2+i\xi_n)c(\xi')c(dx_n)\partial_{x_n}c(\xi')+ic(dx_n)b^2_0c(\xi')
+ic(\xi')b^2_0c(dx_n)-i\partial_{x_n}c(\xi')\big]\nonumber\\
&=&\frac{1}{4(\xi_n-i)^2}\Big[\frac{5}{2}h'(0)c(dx_n)-\frac{5i}{2}h'(0)c(\xi')
  -(2+i\xi_n)c(\xi')c(dx_n)\partial_{\xi_n}c(\xi')+i\partial_{\xi_n}c(\xi')\Big]  ;         \\
B_2&=&\frac{h'(0)}{2}\Big[\frac{c(dx_n)}{4i(\xi_n-i)}+\frac{c(dx_n)-ic(\xi')}{8(\xi_n-i)^2}
+\frac{3\xi_n-7i}{8(\xi_n-i)^3}\big(ic(\xi')-c(dx_n)\big)\Big].
\end{eqnarray}

By (4.30) and (4.40), we have
\begin{eqnarray}
&&{\rm tr }[B_2\times\partial_{\xi_n}\sigma_{-3}((\widehat{D}^{*}\widehat{D}\widehat{D}^{*})^{-1})(x_0)]|_{|\xi'|=1}\nonumber\\
&=&{\rm tr }\Big\{ \frac{h'(0)}{2}\Big[\frac{c(dx_n)}{4i(\xi_n-i)}+\frac{c(dx_n)-ic(\xi')}{8(\xi_n-i)^2}
+\frac{3\xi_n-7i}{8(\xi_n-i)^3}[ic(\xi')-c(dx_n)]\Big] \nonumber\\
&&\times\frac{-4i\xi_nc(\xi')+(i-3i\xi_n^{2})c(dx_n)}{(1+\xi_n^{2})^3}\Big\} \nonumber\\
&=&8h'(0)\frac{4i-11\xi_n-6i\xi_n^{2}+3\xi_n^{3}}{(\xi_n-i)^5(\xi_n+i)^3}.
\end{eqnarray}
Similarly, we have
\begin{eqnarray}
&&{\rm tr }[B_1\times\partial_{\xi_n}\sigma_{-3}((\widehat{D}^{*}\widehat{D}\widehat{D}^{*})^{-1})(x_0)]|_{|\xi'|=1}\nonumber\\
&=&{\rm tr }\Big\{ \frac{1}{4(\xi_n-i)^2}\Big[\frac{5}{2}h'(0)c(dx_n)-\frac{5i}{2}h'(0)c(\xi')
  -(2+i\xi_n)c(\xi')c(dx_n)\partial_{\xi_n}c(\xi')+i\partial_{\xi_n}c(\xi')\Big]\nonumber\\
&&\times \frac{-4i\xi_nc(\xi')+(i-3i\xi_n^{2})c(dx_n)}{(1+\xi_n^{2})^3}\Big\} \nonumber\\
&=&8h'(0)\frac{3+12i\xi_n+3\xi_n^{2}}{(\xi_n-i)^4(\xi_n+i)^3};\\
&&{\rm tr }\bigg[\pi^+_{\xi_n}\Big(\frac{c(\xi)[\bar{c}(\theta)+c(\theta')]
(x_0)c(\xi)}{(1+\xi_n^2)^2}\Big)\times
\partial_{\xi_n}\sigma_{-3}((\widehat{D}^{*}\widehat{D}\widehat{D}^{*})^{-1})
(x_0)\bigg]\bigg|_{|\xi'|=1}\nonumber\\
&=&\frac{2-8i\xi_n-6\xi_n^2}{4(\xi_n-i)^{2}(1+\xi_n^2)^{3}}{\rm tr }[[\bar{c}(\theta)+c(\theta')](x_0)c(\xi')]\nonumber\\
&=&\frac{2-8i\xi_n-6\xi_n^2}{4(\xi_n-i)^{2}(1+\xi_n^2)^{3}}{\rm tr }[c(\theta')(x_0)c(\xi')]\nonumber\\
&=&\frac{2-8i\xi_n-6\xi_n^2}{4(\xi_n-i)^{2}(1+\xi_n^2)^{3}}[-g(\theta',\xi')]
{\rm tr}[\texttt{id}].
\end{eqnarray}

By $\int_{|\xi'|=1}\xi_{1}\cdot\cdot\cdot\xi_{2q+1}\sigma(\xi')=0,$ we have\\
\begin{eqnarray}
{\rm case~(c)}&=&
 -i h'(0)\int_{|\xi'|=1}\int^{+\infty}_{-\infty}
 8\times\frac{-7i+26\xi_n+15i\xi_n^{2}}{(\xi_n-i)^5(\xi_n+i)^3}d\xi_n\sigma(\xi')dx' \nonumber\\
 &&-i\int_{|\xi'|=1}\int^{+\infty}_{-\infty}
 \bigg[\frac{2-8i\xi_n-6\xi_n^2}{4(\xi_n-i)^{2}(1+\xi_n^2)^{3}}[-g(\theta',\xi')]
{\rm tr}[\texttt{id}]\bigg]\bigg|_{|\xi'|=1}d\xi_n\sigma(\xi')dx'\nonumber\\
&=&-8i h'(0)\times\frac{2 \pi i}{4!}\Big[\frac{-7i+26\xi_n+15i\xi_n^{2}}{(\xi_n+i)^3}
     \Big]^{(5)}|_{\xi_n=i}\Omega_4dx'\nonumber\\
&=&\frac{55}{2}\pi h'(0)\Omega_4dx'.
\end{eqnarray}

Now $\Psi$ is the sum of the cases (a), (b) and (c), then
\begin{equation}
\Psi=(\frac{65}{8}-\frac{41}{8}i)\pi h'(0)\Omega_4dx'+120i\pi g(dx_n,\theta')\Omega_4dx'.
\end{equation}

\begin{thm}
Let $M$ be a $6$-dimensional
compact oriented manifold with the boundary $\partial M$ and the metric
$g^M$ as above, $\widehat{D}$ and $\widehat{D}^*$ be modified Novikov operators on $\widehat{M}$, then
\begin{eqnarray}
&&\widetilde{{\rm Wres}}[\pi^+\widehat{D}^{-1}\circ\pi^+(\widehat{D}^{*}\widehat{D}
      \widehat{D}^{*})^{-1}]\nonumber\\
&=&128\pi^3\int_{M}\bigg[64g(\widetilde{e_{j}},
\nabla^{TM}_{\widetilde{e_{j}}}\theta')-\frac{16}{3}s
-64|\theta|^2+256|\theta'|^2\bigg]d{\rm Vol_{M}}\nonumber\\
&&+\int_{\partial M}\bigg[(\frac{65}{8}-\frac{41}{8}i)\pi h'(0)+120i\pi g(dx_n,\theta')\bigg]\Omega_4dx'.
\end{eqnarray}
where $s$ is the scalar curvature.
\end{thm}

On the other hand, we prove the Kastler-Kalau-Walze type theorem for $6$-dimensional manifold with boundary associated to $\widehat{D}^{3}$. An application of (2.1.4) in \cite{Wa5} shows that

\begin{equation}
\widetilde{{\rm Wres}}[\pi^+\widehat{D}^{-1}\circ\pi^+\widehat{D}^{-3}
      ]=\int_M\int_{|\xi|=1}{\rm
trace}_{\wedge^*T^*M}[\sigma_{-4}(\widehat{D}^{-4})]\sigma(\xi)dx+\int_{\partial M}\widehat{\Psi},
\end{equation}
where $\widetilde{{\rm Wres}}$ denote noncommutative residue on minifolds with boundary,
\begin{eqnarray}
\widehat{\Psi} &=&\int_{|\xi'|=1}\int^{+\infty}_{-\infty}\sum^{\infty}_{j, k=0}\sum\frac{(-i)^{|\alpha|+j+k+1}}{\alpha!(j+k+1)!}
\times {\rm trace}_{{\wedge^*T^*M}}[\partial^j_{x_n}\partial^\alpha_{\xi'}\partial^k_{\xi_n}\sigma^+_{r}(\widehat{D}^{-1})(x',0,\xi',\xi_n)
\nonumber\\
&&\times\partial^\alpha_{x'}\partial^{j+1}_{\xi_n}\partial^k_{x_n}\sigma_{l}
(\widehat{D}^{-3})(x',0,\xi',\xi_n)]d\xi_n\sigma(\xi')dx',
\end{eqnarray}
and the sum is taken over $r+\ell-k-j-|\alpha|-1=-6, \ r\leq-1, \ell\leq -3$.

Locally we can use Theorem 2.2 (2.29) to compute the interior term of (4.48), we have
\begin{eqnarray}
&&\int_M\int_{|\xi|=1}{\rm
trace}_{\wedge^*T^*M}[\sigma_{-4}(\widehat{D}^{-4})]\sigma(\xi)dx
=128\pi^3\int_{M}
\bigg[-\frac{16}{3}s
-64|\theta|^2\bigg]d{\rm Vol_{M}}.
\end{eqnarray}

So we only need to compute $\int_{\partial M} \widehat{\Psi}$. Let us now turn to compute the specification of
$\widehat{D}^3$.
\begin{eqnarray}
\widehat{D}^3
&=&\sum^{n}_{i=1}c(\widetilde{e_{i}})\langle \widetilde{e_{i}},dx_{l}\rangle(-g^{ij}\partial_{l}\partial_{i}\partial_{j})
+\sum^{n}_{i=1}c(\widetilde{e_{i}})\langle \widetilde{e_{i}},dx_{l}\rangle \bigg\{-(\partial_{l}g^{ij})\partial_{i}\partial_{j}-g^{ij}\bigg(4(\sigma_{i}
+a_{i})\partial_{j}-2\Gamma^{k}_{ij}\partial_{k}\bigg)\partial_{l}\bigg\} \nonumber\\
&&+\sum^{n}_{i=1}c(\widetilde{e_{i}})\langle \widetilde{e_{i}},dx_{l}\rangle \bigg\{-2(\partial_{l}g^{ij})(\sigma_{i}+a_i)\partial_{j}+g^{ij}
(\partial_{l}\Gamma^{k}_{ij})\partial_{k}-2g^{ij}[(\partial_{l}\sigma_{i})
+(\partial_{l}a_i)]\partial_{j}
         +(\partial_{l}g^{ij})\Gamma^{k}_{ij}\partial_{k}\nonumber\\
         &&+\sum_{j,k}\Big[\partial_{l}\Big(c(\theta')c(\widetilde{e_{j}})
         +c(\widetilde{e_{j}}) c(\theta')\Big)\Big]\langle \widetilde{e_{j}},dx^{k}\rangle\partial_{k}
         +\sum_{j,k}\Big(c(\theta')c(\widetilde{e_{j}})+c(\widetilde{e_{j}})c(\theta')\Big)\Big[\partial_{l}\langle \widetilde{e_{j}},dx^{k}\rangle\Big]\partial_{k} \bigg\}\nonumber\\
         &&+\sum^{n}_{i=1}c(\widetilde{e_{i}})\langle \widetilde{e_{i}},dx_{l}\rangle\partial_{l}\bigg\{-g^{ij}\Big[
         (\partial_{i}\sigma_{j})+(\partial_{i}a_{j})+\sigma_{i}\sigma_{j}+\sigma_{i}a_{j}
         +a_{i}\sigma_{j}+a_{i}a_{j}-\Gamma_{i,j}^{k}\sigma_{k}-\Gamma_{i,j}^{k}a_{k}\nonumber\\
         &&+\sum_{i,j}g^{i,j}\Big[c(\theta')c(\partial_{i})\sigma_{i}
         +c(\theta')c(\partial_{i})a_{i}+c(\partial_{i})\partial_{i}(c(\theta'))
         +c(\partial_{i})\sigma_{i}c(\theta')+c(\partial_{i})a_{i}c(\theta')\Big]+\frac{1}{4}s\nonumber\\
         &&-\frac{1}{8}\sum_{ijkl}R_{ijkl}\bar{c}(\widetilde{e_i})\bar{c}(\widetilde{e_j})
         c(\widetilde{e_k})c(\widetilde{e_l})
         +\sum_{i}c(\widetilde{e_i})\bar{c}(\nabla_{\widetilde{e_i}}^{TM}\theta)
         +|\theta|^{2}-|\theta'|^2+\bar{c}(\theta)c(\theta')+c(\theta')\bar{c}(\theta)\bigg\}\nonumber\\
         &&+\Big[(\sigma_{i}+a_{i})+(\bar{c}(\theta)+c(\theta'))\Big](-g^{ij}\partial_{i}\partial_{j})
         +\sum^{n}_{i=1}c(\widetilde{e_{i}})\langle \widetilde{e_{i}},dx_{l}\rangle \bigg\{2\sum_{j,k}\Big[c(\theta')c(\widetilde{e_{j}})+c(\widetilde{e_{j}})c(\theta')\Big]\nonumber\\
         &&\times\langle \widetilde{e_{i}},dx_{k}\rangle\bigg\}\partial_{l}\partial_{k}
         +\Big[(\sigma_{i}+a_{i})+(\bar{c}(\theta)+c(\theta'))\Big]
         \bigg\{-\sum_{i,j}g^{i,j}\Big[2\sigma_{i}\partial_{j}+2a_{i}\partial_{j}
         -\Gamma_{i,j}^{k}\partial_{k}+(\partial_{i}\sigma_{j})\nonumber\\
         &&+(\partial_{i}a_{j})+\sigma_{i}\sigma_{j}+\sigma_{i}a_{j}+a_{i}\sigma_{j}+a_{i}a_{j} -\Gamma_{i,j}^{k}\sigma_{k}-\Gamma_{i,j}^{k}a_{k}\Big]+\sum_{i,j}g^{i,j}\Big[c(\partial_{i})c(\theta')
         +c(\theta')c(\partial_{i})\Big]\partial_{j}\nonumber\\
         &&+\sum_{i,j}g^{i,j}\Big[c(\theta')c(\partial_{i})\sigma_{i}+c(\theta')c(\partial_{i})a_{i}
         +c(\partial_{i})\partial_{i}(c(\theta'))+c(\partial_{i})\sigma_{i}c(\theta')
         +c(\partial_{i})a_{i}c(\theta')\Big]+\frac{1}{4}s-|\theta'|^2\nonumber\\
         &&-\frac{1}{8}\sum_{ijkl}R_{ijkl}\bar{c}(\widetilde{e_i})\bar{c}(\widetilde{e_j})
         c(\widetilde{e_k})c(\widetilde{e_l})
         +\sum_{i}c(\widetilde{e_i})\bar{c}(\nabla_{\widetilde{e_i}}^{TM}\theta)
         +|\theta|^{2}+\bar{c}(\theta)c(\theta')+c(\theta')\bar{c}(\theta)\bigg\}.
\end{eqnarray}

Then, we obtain
\begin{lem} The following identities hold:
\begin{eqnarray}
\sigma_2(\widehat{D}^3)&=&\sum_{i,j,l}c(dx_{l})\partial_{l}(g^{i,j})\xi_{i}\xi_{j} +c(\xi)(4\sigma^k+4a^k-2\Gamma^k)\xi_{k}-2[c(\xi)c(\theta')c(\xi)-|\xi|^2c(\theta')]\nonumber\\
&&+\frac{1}{4}|\xi|^2\sum_{s,t,l}\omega_{s,t}
(\widetilde{e_l})[c(\widetilde{e_l})\bar{c}(\widetilde{e_s})\bar{c}(\widetilde{e_t})
-c(\widetilde{e_l})c(\widetilde{e_s})c(\widetilde{e_t})]
+|\xi|^2(\bar{c}(\theta)+c(\theta'));\nonumber\\
\sigma_{3}(\widehat{D}^3)&=&ic(\xi)|\xi|^{2}.
\end{eqnarray}
\end{lem}

Write
\begin{eqnarray}
\sigma(\widehat{D}^3)&=&p_3+p_2+p_1+p_0;
~\sigma(\widehat{D}^{-3})=\sum^{\infty}_{j=3}q_{-j}.
\end{eqnarray}

By the composition formula of pseudodifferential operators, we have

\begin{eqnarray}
1=\sigma(\widehat{D}^3\circ \widehat{D}^{-3})&=&
\sum_{\alpha}\frac{1}{\alpha!}\partial^{\alpha}_{\xi}
[\sigma(\widehat{D}^3)]D^{\alpha}_{x}
[\sigma(\widehat{D}^{-3})] \nonumber\\
&=&(p_3+p_2+p_1+p_0)(q_{-3}+q_{-4}+q_{-5}+\cdots) \nonumber\\
&+&\sum_j(\partial_{\xi_j}p_3+\partial_{\xi_j}p_2+\partial_{\xi_j}p_1+\partial_{\xi_j}p_0)
(D_{x_j}q_{-3}+D_{x_j}q_{-4}+D_{x_j}q_{-5}+\cdots) \nonumber\\
&=&p_3q_{-3}+(p_3q_{-4}+p_2q_{-3}+\sum_j\partial_{\xi_j}p_3D_{x_j}q_{-3})+\cdots,
\end{eqnarray}
by (4.53), we have

\begin{equation}
q_{-3}=p_3^{-1};~q_{-4}=-p_3^{-1}[p_2p_3^{-1}+\sum_j\partial_{\xi_j}p_3D_{x_j}(p_3^{-1})].
\end{equation}
By (4.50)-(4.54), we have some symbols of operators.
\begin{lem} The following identities hold:
\begin{eqnarray}
\sigma_{-3}(\widehat{D}^{-3})&=&\frac{ic(\xi)}{|\xi|^{4}};\nonumber\\
\sigma_{-4}(\widehat{D}^{-3})&=&
\frac{c(\xi)\sigma_{2}(\widehat{D}^{3})c(\xi)}{|\xi|^8}
+\frac{ic(\xi)}{|\xi|^8}\Big(|\xi|^4c(dx_n)\partial_{x_n}c(\xi')
-2h'(0)c(dx_n)c(\xi)\nonumber\\
&&+2\xi_{n}c(\xi)\partial_{x_n}c(\xi')+4\xi_{n}h'(0)\Big).
\end{eqnarray}
\end{lem}

From the remark above, now we can compute $\widehat{\Psi}$ (see formula (4.48) for the definition of $\widehat{\Psi}$). We use ${\rm tr}$ as shorthand of ${\rm trace}$. Since $n=6$, then ${\rm tr}_{\wedge^*T^*M}[\texttt{id}]=64$.
Since the sum is taken over $r+\ell-k-j-|\alpha|-1=-6, \ r\leq-1, \ell\leq -3$, then we have the $\int_{\partial_{M}}\widehat{\Psi}$
is the sum of the following five cases:

~\\
\noindent  {\bf case (a)~(I)}~$r=-1, l=-3, j=k=0, |\alpha|=1$.\\
By (4.48), we get
 \begin{equation}
{\rm case~(a)~(I)}=-\int_{|\xi'|=1}\int^{+\infty}_{-\infty}\sum_{|\alpha|=1}{\rm trace}
\Big[\partial^{\alpha}_{\xi'}\pi^{+}_{\xi_{n}}\sigma_{-1}(\widehat{D}^{-1})
      \times\partial^{\alpha}_{x'}\partial_{\xi_{n}}\sigma_{-3}(\widehat{D}^{-3})\Big](x_0)d\xi_n\sigma(\xi')dx'.
\end{equation}

\noindent  {\bf case (a)~(II)}~$r=-1, l=-3, |\alpha|=k=0, j=1$.\\
By (4.48), we have
  \begin{equation}
{\rm case~(a)~(II)}=-\frac{1}{2}\int_{|\xi'|=1}\int^{+\infty}_{-\infty} {\rm
trace} \Big[\partial_{x_{n}}\pi^{+}_{\xi_{n}}\sigma_{-1}(\widehat{D}^{-1})
      \times\partial^{2}_{\xi_{n}}\sigma_{-3}(\widehat{D}^{-3})\Big](x_0)d\xi_n\sigma(\xi')dx'.
\end{equation}

\noindent  {\bf case (a)~(III)}~$r=-1,l=-3,|\alpha|=j=0,k=1$.\\
By (4.48), we have
 \begin{equation}
{\rm case~ (a)~(III)}=-\frac{1}{2}\int_{|\xi'|=1}\int^{+\infty}_{-\infty}{\rm trace} \Big[\partial_{\xi_{n}}\pi^{+}_{\xi_{n}}\sigma_{-1}(\widehat{D}^{-1})
      \times\partial_{\xi_{n}}\partial_{x_{n}}\sigma_{-3}(\widehat{D}^{-3})\Big](x_0)d\xi_n\sigma(\xi')dx'.
\end{equation}\\
By Lemma 4.2 and Lemma 4.5, we have $\sigma_{-3}((\widehat{D}^*\widehat{D}\widehat{D}^*)^{-1})=\sigma_{-3}(\widehat{D}^{-3})$
, by (4.10)-(4.20), we obtain
$${\bf case~(a)}=5\pi h'(0)\Omega_{4}dx',$$
 where ${\rm \Omega_{4}}$ is the canonical volume of $S^{4}.$\\

\noindent  {\bf case (b)}~$r=-1,l=-4,|\alpha|=j=k=0$.\\
By (4.48), we have
 \begin{eqnarray}
{\rm case~ (b)}&=&-i\int_{|\xi'|=1}\int^{+\infty}_{-\infty}{\rm trace} \Big[\pi^{+}_{\xi_{n}}\sigma_{-1}(\widehat{D}^{-1})
      \times\partial_{\xi_{n}}\sigma_{-4}(\widehat{D}^{-3})\Big](x_0)d\xi_n\sigma(\xi')dx'\nonumber\\
&=&i\int_{|\xi'|=1}\int^{+\infty}_{-\infty}{\rm trace} [\partial_{\xi_n}\pi^+_{\xi_n}\sigma_{-1}(\widehat{D}^{-1})\times
\sigma_{-4}(\widehat{D}^{-3})](x_0)d\xi_n\sigma(\xi')dx'.
\end{eqnarray}

In the normal coordinate, $g^{ij}(x_{0})=\delta^{j}_{i}$ and $\partial_{x_{j}}(g^{\alpha\beta})(x_{0})=0$, if $j<n$; $\partial_{x_{j}}(g^{\alpha\beta})(x_{0})=h'(0)\delta^{\alpha}_{\beta}$, if $j=n$.
So by Lemma A.2 in \cite{Wa3}, we have $\Gamma^{n}(x_{0})=\frac{5}{2}h'(0)$ and $\Gamma^{k}(x_{0})=0$ for $k<n$. By the definition of $\delta^{k}$ and Lemma 2.3 in \cite{Wa3}, we have $\delta^{n}(x_{0})=0$ and $\delta^{k}=\frac{1}{4}h'(0)c(\widetilde{e_{k}})c(\widetilde{e_{n}})$ for $k<n$. By Lemma 4.5, we obtain

\begin{eqnarray}
\sigma_{-4}(\widehat{D}^{-3})(x_{0})|_{|\xi'|=1}&=&
\frac{c(\xi)\sigma_{2}(\widehat{D}^{-3})
(x_{0})|_{|\xi'|=1}c(\xi)}{|\xi|^8}
-\frac{c(\xi)}{|\xi|^4}\sum_j\partial_{\xi_j}\big(c(\xi)|\xi|^2\big)
D_{x_j}\big(\frac{ic(\xi)}{|\xi|^4}\big)\nonumber\\
&=&\frac{1}{|\xi|^8}c(\xi)\Big(\frac{1}{2}h'(0)c(\xi)\sum_{k<n}\xi_k
c(\widetilde{e_k})c(\widetilde{e_n})-\frac{1}{2}h'(0)c(\xi)\sum_{k<n}\xi_k
\bar{c}(\widetilde{e_k})\bar{c}(\widetilde{e_n})\nonumber\\
&&-\frac{5}{2}h'(0)\xi_nc(\xi)-\frac{1}{4}h'(0)|\xi|^2c(dx_n)
-2[c(\xi)c(\theta')c(\xi)-|\xi|^2c(\theta')]\nonumber\\
&&+|\xi|^2(\widehat{c}(\theta)+c(\theta'))\Big)c(\xi)
+\frac{ic(\xi)}{|\xi|^8}\Big(|\xi|^4c(dx_n)\partial_{x_n}c(\xi')
-2h'(0)c(dx_n)c(\xi)\nonumber\\
&&+2\xi_{n}c(\xi)\partial_{x_n}c(\xi')+4\xi_{n}h'(0)\Big).
\end{eqnarray}
By (3.34) and (4.60), we have
\begin{eqnarray}
&&{\rm tr} [\partial_{\xi_n}\pi^+_{\xi_n}\sigma_{-1}(\widehat{D}^{-1})\times
\sigma_{-4}(\widehat{D}^{-3}) ](x_0)|_{|\xi'|=1} \nonumber\\
%&&={\rm tr} \Big[-\frac{c(\xi')+ic(dx_n)}{2(\xi_n-i)^2} \times
 %\frac{1}{(1+\xi_{n}^{2})^{4}}\Big\{\Big( (3+4i)\xi_{n}+3\xi_{n}^{3}\Big)h'(0)c(\xi') \nonumber\\
 %&& +\big(\frac{3}{4}-2i+(6i+2)\xi_{n}^{2}+\frac{9}{4}\xi_{n}^{4}\big)h'(0)c(dx_n) \nonumber\\
%&& +\big(-\frac{1}{2}\xi_{n}-\frac{1}{2}\xi_{n}^{3}\big)h'(0)\bar{c}(\xi')\bar{c}(dx_n)c(dx_n) \nonumber\\
%&& +\big(-\frac{1}{2}-\frac{1}{2}\xi_{n}^{2}\big)h'(0)\bar{c}(\xi')\bar{c}(dx_n)c(\xi') \nonumber\\
%&& +\big(1+2\xi_{n}^{2}+\xi_{n}^{4}\big)T\bar{c}(v)+\big(-3i\xi_{n}-4i\xi_{n}^{3}-i\xi_{n}^{5}\big)\partial_{x_n}c(\xi')\big\}\nonumber\\
&=&\frac{1}{2(\xi_{n}-i)^{2}(1+\xi_{n}^{2})^{4}}\big(\frac{3}{4}i+2+(3+4i)\xi_{n}+(-6+2i)\xi_{n}^{2}+3\xi_{n}^{3}+\frac{9i}{4}\xi_{n}^{4}\big)h'(0){\rm tr}
[id]\nonumber\\
&&+\frac{1}{2(\xi_{n}-i)^{2}(1+\xi_{n}^{2})^{4}}\big(-1-3i\xi_{n}-2\xi_{n}^{2}-4i\xi_{n}^{3}-\xi_{n}^{4}-i\xi_{n}^{5}\big){\rm tr[c(\xi')\partial_{x_n}c(\xi')]}\nonumber\\
&&-\frac{1}{2(\xi_{n}-i)^{2}(1+\xi_{n}^{2})^{4}}\big(\frac{1}{2}i+\frac{1}{2}\xi_{n}+\frac{1}{2}\xi_{n}^{2}+\frac{1}{2}\xi_{n}^{3}\big){\rm tr}
[c(\xi')\bar{c}(\xi')c(dx_n)\bar{c}(dx_n)]\nonumber\\
&&+\frac{-3\xi_ni+1}{2(\xi_{n}-i)^{4}(i+\xi_{n})^{3}}{\rm tr}\big[c(\theta')c(dx_n)\big]
-\frac{\xi_n+3i}{2(\xi_{n}-i)^{4}(i+\xi_{n})^{3}}{\rm tr}\big[c(\theta')c(\xi')\big].\nonumber\\
\end{eqnarray}
By (4.24),(4.25), we have
\begin{eqnarray}
{\rm case~ (b)}&=&
 ih'(0)\int_{|\xi'|=1}\int^{+\infty}_{-\infty}64\times\frac{\frac{3}{4}i+2+(3+4i)\xi_{n}+(-6+2i)\xi_{n}^{2}+3\xi_{n}^{3}+\frac{9i}{4}\xi_{n}^{4}}{2(\xi_n-i)^5(\xi_n+i)^4}d\xi_n\sigma(\xi')dx'\nonumber\\ &+&ih'(0)\int_{|\xi'|=1}\int^{+\infty}_{-\infty}32\times\frac{1+3i\xi_{n}
 +2\xi_{n}^{2}+4i\xi_{n}^{3}+\xi_{n}^{4}+i\xi_{n}^{5}}{2(\xi_{n}-i)^{2}
 (1+\xi_{n}^{2})^{4}}d\xi_n\sigma(\xi')dx'\nonumber\\
 &+&i\int_{|\xi'|=1}\int^{+\infty}_{-\infty}
 \frac{-3\xi_ni +1}{2(\xi_{n}-i)^{4}(i+\xi_{n})^{3}}{\rm tr}\big[c(\theta')c(dx_n)\big]d\xi_n\sigma(\xi')dx'\nonumber\\
 &-&i\int_{|\xi'|=1}\int^{+\infty}_{-\infty}
 \frac{\xi_n+3i}{2(\xi_{n}-i)^{4}(i+\xi_{n})^{3}}{\rm tr}\big[c(\theta')c(\xi')\big]d\xi_n\sigma(\xi')dx'\nonumber\\
&=&(-\frac{41}{8}i-\frac{195}{8})\pi h'(0)\Omega_4dx'.
\end{eqnarray}

\noindent {\bf  case (c)}~$r=-2,l=-3,|\alpha|=j=k=0$.\\
By (4.48), we have

\begin{equation}
{\rm case~ (c)}=-i\int_{|\xi'|=1}\int^{+\infty}_{-\infty}{\rm trace} \Big[\pi^{+}_{\xi_{n}}\sigma_{-2}(\widehat{D}^{-1})
      \times\partial_{\xi_{n}}\sigma_{-3}(\widehat{D}^{-3})\Big](x_0)d\xi_n\sigma(\xi')dx'.
\end{equation}\\
By Lemma 4.2 and Lemma 4.5, we have $\sigma_{-3}((\widehat{D}^*\widehat{D}\widehat{D}^*)^{-1})=\sigma_{-3}(\widehat{D}^{-3})$
, by (4.27)-(4.44), we obtain

$$
{\rm case~ (c)}=\frac{55}{2}\pi h'(0)\Omega_4dx'.
$$

Now $\widehat{\Psi}$ is the sum of the cases (a), (b) and (c), then
\begin{equation}
\widehat{\Psi}=\bigg[(\frac{65}{8}-\frac{41}{8}i)\pi h'(0)\bigg]\Omega_4dx'.
\end{equation}

\begin{thm}
Let $M$ be a $6$-dimensional
compact oriented manifold with the boundary $\partial M$ and the metric
$g^M$ as above, $\widehat{D}$ be a modified Novikov operator on $\widehat{M}$, then
\begin{eqnarray}
&&\widetilde{{\rm Wres}}[\pi^+\widehat{D}^{-1}\circ\pi^+(
      \widehat{D}^{-3})]\nonumber\\
&=&128\pi^3\int_{M}
\bigg[-\frac{16}{3}s
-64|\theta|^2\bigg]d{\rm Vol_{M}}+\int_{\partial M}\bigg[(\frac{65}{8}-\frac{41}{8}i)\pi h'(0)\bigg]\Omega_4dx'.
\end{eqnarray}
where $s$ is the scalar curvature.
\end{thm}

%%%%%%%%%第5部分%%%%%%%
\section{The Spectral Action for Witten Deformation}

In this section, we will compute the spectral action for Witten deformation. Let $(M, g^{M})$ be a $n$-dimensional oriented compact Riemannian manifold. Now we will recall the definition of the Witten deformation $D_{\theta}$ (see details in \cite{wpz}).

Let $\nabla^L$ denote the Levi-Civita connection about $g^M$ which is Riemannian metric of $M$. In the local coordinates $\{x_i; 1\leq i\leq n\}$ and the
fixed orthonormal frame $\{\widetilde{e_1},\cdots,\widetilde{e_n}\}$, the connection matrix $(\omega_{s,t})$ is defined by
\begin{equation}
\nabla^L(\widetilde{e_1},\cdots,\widetilde{e_n})= (\widetilde{e_1},\cdots,\widetilde{e_n})(\omega_{s,t}).
\end{equation}
 Let $\epsilon (\widetilde{e_j*} ),~\iota (\widetilde{e_j*} )$ be the exterior and interior multiplications respectively. Write
\begin{equation}
c(\widetilde{e_j})=\epsilon (\widetilde{e_j*} )-\iota (\widetilde{e_j*} );~~
\bar{c}(\widetilde{e_j})=\epsilon (\widetilde{e_j*} )+\iota
(\widetilde{e_j*} ).
\end{equation}
 The Witten deformation is defined by
\begin{equation}
D_{\theta}=d+\delta+\bar{c}(\theta)=\sum^n_{i=1}c(\widetilde{e_i})\bigg[\widetilde{e_i}
+\frac{1}{4}\sum_{s,t}\omega_{s,t}
(\widetilde{e_i})[\bar{c}(\widetilde{e_s})\bar{c}(\widetilde{e_t})-c(\widetilde{e_s})
c(\widetilde{e_t})]\bigg]+\bar{c}(\theta).
\end{equation}
where $d,~\delta$, $\theta\in\Gamma(TM)$.

By Proposition 4.6 of \cite{wpz}, we have
\begin{equation}
D^{2}_{\theta}=(d+\delta)^{2}+\sum_{i}c(\widetilde{e_i})\bar{c}(
\nabla_{\widetilde{e_i}}^{TM}\theta)
+|\theta|^{2}.
\end{equation}

By \cite{Y}, the local expression of $(d+\delta)^{2}$ is
\begin{equation}
(d+\delta)^{2}
=-\triangle_{0}-\frac{1}{8}\sum_{ijkl}R_{ijkl}\bar{c}(\widetilde{e_i})\bar{c}(\widetilde{e_j})c(\widetilde{e_k})c(\widetilde{e_l})+\frac{1}{4}s.
\end{equation}

Let $g^{ij}=g(dx_{i},dx_{j})$, $\xi=\sum_{k}\xi_{j}dx_{j}$ and $\nabla^L_{\partial_{i}}\partial_{j}=\sum_{k}\Gamma_{ij}^{k}\partial_{k}$,  we denote
\begin{eqnarray}
&&\sigma_{i}=-\frac{1}{4}\sum_{s,t}\omega_{s,t}
(\widetilde{e_i})c(\widetilde{e_s})c(\widetilde{e_t})
;~~~a_{i}=\frac{1}{4}\sum_{s,t}\omega_{s,t}
(\widetilde{e_i})\bar{c}(\widetilde{e_s})\bar{c}(\widetilde{e_t});\nonumber\\
&&\xi^{j}=g^{ij}\xi_{i},~~~~\Gamma^{k}=g^{ij}\Gamma_{ij}^{k},~~~~\sigma^{j}=g^{ij}\sigma_{i},~~~~a^{j}=g^{ij}a_{i}.
\end{eqnarray}

For a smooth vector field $X$ on $M$, let $c(X)$
denote the Clifford action. Since $E$ is globally
defined on $M$, so we can preform computations of $E$ in normal coordinates. Taking normal coordinates about $x_0$, then $\sigma^{i}(x_0)=0$, $a^{i}(x_0)=0$, $\partial^{j}[c(\partial_{j})](x_0)=0$,
$\Gamma^k(x_0)=0$, $g^{ij}(x_0)=\delta^j_i$, so that
\begin{eqnarray}
E(x_0)&=&\frac{1}{8}\sum_{ijkl}R_{ijkl}\bar{c}(\widetilde{e_i})\bar{c}
(\widetilde{e_j})c(\widetilde{e_k})c(\widetilde{e_l})
-\frac{1}{4}s-\sum_{i}c(\widetilde{e_i})\bar{c}
(\nabla_{\widetilde{e_i}}^{TM}\theta)-|\theta|^{2}.
\end{eqnarray}
For the Witten deformation $D_{\theta}$, we will compute the spectral action for it on 4-dimensional compact manifold. We will calculate the bosonic part of the spectral action for the Witten deformation.  It is defined to be the number of eigenvalues of $D_{\theta}$ in the interval $[-\wedge,\wedge]$ with $\wedge\in \mathbb{R}^+$. It is expressed as
\begin{equation}
I={\rm tr}\widehat{F}(\frac{D^2_{\theta}}{\wedge^2}),
\end{equation}
Here, ${\rm tr}$ denotes the operator trace in the $L^2$ completion of $\Gamma(M,S(TM))$ and $\widehat{F}:\mathbb{R}^+\rightarrow\mathbb{R}^+$ is a cut-off function with support in the interval $[0,1]$ which is constant near the origin. By Lemma 1.7.4 in \cite{G}, we have the heat trace asymptotics, for $t\rightarrow0$,
\begin{equation}
{\rm tr}(e^{-tD^2_{\theta}})\thicksim\sum_{m\geq0}t^{m-\frac{n}{2}}
a_{2m}(D^2_{\theta}),
\end{equation}
One uses the Seely-DeWitt coefficients $a_{2m}(D^2_{\theta})$ and $t=\wedge^{-2}$ to obtain an asympotics for the spectral action when ${\rm {dim}}M=4$,
\begin{equation}
I={\rm tr}\widehat{F}(\frac{D^2_{\theta}}{\wedge^2})
\sim\wedge^{4}F_4a_0(D^2_{\theta})+\wedge^{2}F_2a_2(D^2_{\theta})+
\wedge^{0}F_0a_4(D^2_{\theta})~~as~~\wedge\rightarrow\infty,
\end{equation}
with the first three moments of the cut-off function which are given by $F_4=\int^\infty_{0}s\widehat{F}(s)ds$, $F_2=\int^\infty_{0}\widehat{F}(s)ds$, $F_0=\widehat{F}(0)$.

We use Theorem 4.1.6 in [22] to obtain the first
three coefficients of the heat trace asymptotics:
\begin{eqnarray}
a_0(D^2_{\theta})
&=&(4\pi)^{-\frac{4}{2}}\int_M{\rm tr}(\texttt{id})dvol\nonumber\\
a_2(D^2_{\theta})
&=&(4\pi)^{-\frac{4}{2}}\int_M{\rm tr}[\frac{s}{6}+E]dvol\nonumber\\
a_4(D^2_{\theta})
&=&\frac{(4\pi)^{-\frac{4}{2}}}
{360}\int_M{\rm tr}[-12R_{ijij,kk}+5R_{ijij}R_{klkl}-2R_{ijik}
R_{ljlk}+2R_{ijkl}R_{ijkl}-60R_{ijij}E\nonumber\\
&&+180E^2+60E_{,kk}+30\Omega_{ij}\Omega_{ij}]dvol,
\end{eqnarray}
By Clifford action and cyclicity of the trace, we have
\begin{eqnarray}
{\rm tr}(c(e_i))=0;{\rm tr}(c(e_i)c(e_j))=0~(i\neq j);
{\rm tr}[\bar{c}(\widetilde{e_i})\bar{c}
(\widetilde{e_j})c(\widetilde{e_k})c(\widetilde{e_l})]=0~(i\neq j);
{\rm tr}[\sum_{i}c(\widetilde{e_i})\bar{c}
(\nabla_{\widetilde{e_i}}^{TM}\theta)]=0.
\end{eqnarray}
So we obtain
\begin{eqnarray}
a_0(D^2_{\theta})
&=&(4\pi)^{-\frac{4}{2}}\int_M{\rm tr}(\texttt{id})dvol=\pi^{-2}\int_Mdvol,\nonumber\\
a_2(D^2_{\theta})
&=&(4\pi)^{-\frac{4}{2}}\int_M{\rm tr}[\frac{s}{6}+E]dvol\nonumber\\
&=&(4\pi)^{-\frac{4}{2}}\int_M{\rm tr}[\frac{s}{6}+\frac{1}{8}\sum_{ijkl}R_{ijkl}\bar{c}(\widetilde{e_i})\bar{c}
(\widetilde{e_j})c(\widetilde{e_k})c(\widetilde{e_l})
-\frac{1}{4}s-\sum_{i}c(\widetilde{e_i})\bar{c}
(\nabla_{\widetilde{e_i}}^{TM}\theta)-|\theta|^{2}]dvol\nonumber\\
&=&-\int_M(\frac{s}{12\pi^{2}}+\frac{|\theta|^{2}}{\pi^{2}})dvol.
\end{eqnarray}
And we have
\begin{eqnarray}
&&\int_M{\rm tr}(5R_{ijij}R_{klkl}-60R_{ijij}E+180E^2-12R_{ijij,kk}+60E_{,kk})dvol\nonumber\\
&=&\int_M\bigg[{\rm tr}(5s^2+60sE+180E^2)-{\rm tr}[12\Delta s]+60[-\Delta({\rm tr}E)]\bigg]dvol\nonumber\\
&=&\int_M{\rm trace}\bigg[\frac{5}{4}s^2+180\bigg(\sum_{i}c(\widetilde{e_i})\bar{c}
(\nabla_{\widetilde{e_i}}^{TM}\theta)\bigg)^2+\frac{45}{16}
\bigg(\sum_{ijkl}R_{ijkl}\bar{c}(\widetilde{e_i})\bar{c}
(\widetilde{e_j})c(\widetilde{e_k})c(\widetilde{e_l})\bigg)^2\nonumber\\
&&+30s|\theta|^{2}+180|\theta|^{4}
-45\sum_{ijkl}R_{ijkl}
\bar{c}(\widetilde{e_i})\bar{c}
(\widetilde{e_j})c(\widetilde{e_k})
c(\widetilde{e_l})\sum_{p}
c(\widetilde{e_p})\bar{c}
(\nabla_{\widetilde{e_p}}^{TM}\theta)+48\Delta s+960\Delta(|\theta|^{2})\bigg]dvol\nonumber\\
&=&\int_M\bigg[20s^2+2880(\sum_{i}|\nabla_{\widetilde{e_i}}^{TM}
\theta|^2+|\theta|^{4})
+180\sum_{ijkl}R^2_{ijkl}+480s|\theta|^{2}\bigg]dvol.
\end{eqnarray}

And
${\rm tr}[\Omega_{ij}\Omega_{ij}]$ is globally defined, so we only compute it in normal coordinates about $x_0$ and the local orthonormal frame $e_i$ obtained by parallel transport along geodesics from $x_0$. Then
\begin{eqnarray}
\omega_{s,t}(x_0)=0,~~\partial_i(c(\widetilde{e_j}))=0,
~~[\widetilde{e_i},\widetilde{e_j}](x_0)=0.
\end{eqnarray}

Then, we have
\begin{eqnarray}
\Omega(\widetilde{e_i},\widetilde{e_j})(x_0)
&=&\nabla^{\wedge^*T^*M}_{\widetilde{e}_i}\nabla^{\wedge^*T^*M}_{\widetilde{e}_j}
-\nabla^{\wedge^*T^*M}_{\widetilde{e}_j}\nabla^{\wedge^*T^*M}_{\widetilde{e}_i}
-\nabla^{\wedge^*T^*M}_{[\widetilde{e}_i,\widetilde{e}_j]}\nonumber\\
&=&-\frac{1}{4}\sum^{n}_{s,t=1}R^{M}_{ijst}
[\bar{c}(\widetilde{e_s})\bar{c}(\widetilde{e_t})
-c(\widetilde{e_s})c(\widetilde{e_t})].
\end{eqnarray}

So we have
\begin{eqnarray}
{\rm tr}[\Omega_{ij}\Omega_{ij}](x_0)
&=&{\rm tr}\bigg[\frac{1}{16}\sum_{s,t,s_1,t_1=1}
R^{M}_{ijst}R^{M}_{ijs_1t_1}[\bar{c}(\widetilde{e_s})\bar{c}(\widetilde{e_t})
-c(\widetilde{e_s})c(\widetilde{e_t})][\bar{c}(\widetilde{e_{s_1}})\bar{c}(\widetilde{e_{t_1}})
-c(\widetilde{e_{s_1}})c(\widetilde{e_{t_1}})]\bigg]\nonumber\\
&=&{\rm tr}\bigg[\frac{1}{16}\sum_{s,t,s_1,t_1=1}R^{M}_{ijst}R^{M}_{ijs_1t_1}
[\bar{c}(\widetilde{e_s})\bar{c}(\widetilde{e_t})\bar{c}(\widetilde{e_{s_1}})\bar{c}(\widetilde{e_{t_1}})
-\bar{c}(\widetilde{e_s})\bar{c}(\widetilde{e_t})c(\widetilde{e_{s_1}})c(\widetilde{e_{t_1}})\nonumber\\
&&-c(\widetilde{e_s})c(\widetilde{e_t})\bar{c}(\widetilde{e_{s_1}})\bar{c}(\widetilde{e_{t_1}})
+c(\widetilde{e_s})c(\widetilde{e_t})c(\widetilde{e_{s_1}})c(\widetilde{e_{t_1}})\bigg]\nonumber\\
&=&-4\sum_{ijst}(R^{M}_{ijst})^2
\end{eqnarray}

\begin{prop}
The following equality holds:
\begin{eqnarray}
a_4(D^2_{\theta})
&=&\frac{1}
{5760\pi^2}\int_M[20s^2+2880(\sum_{i}|\nabla_{\widetilde{e_i}}^{TM}
\theta|^2+|\theta|^{4})
+180\sum_{ijkl}R^2_{ijkl}+480s|\theta|^{2}\nonumber\\
&&-32R_{ijik}R_{ljlk}
+32R^2_{ijik}-1920\sum_{i,j,s,t}(R^{M}_{ijst})^2\bigg]dvol
\end{eqnarray}
where $s$ is the scalar curvature.
\end{prop}

\section*{Acknowledgements}
This work was supported by NSFC. 11771070 .
 The authors thank the referee for his (or her) careful reading and helpful comments.

\end{document}